\documentclass[12pt,fleqn]{article}
\usepackage{latexsym}
\usepackage{amsmath}
\usepackage{amssymb}
\usepackage{epsfig}
\usepackage{amsfonts}

\begin{document}

\newtheorem{definition}{Definition}[section]
\newtheorem{theorem}[definition]{Theorem}
\newtheorem{lemma}[definition]{Lemma}
\newtheorem{remark}[definition]{Remark}
\newtheorem{remarks}[definition]{Remarks}
\newtheorem{leeg}[definition]{}
\newtheorem{fact}[definition]{Fact}
\newtheorem{definitions}[definition]{Definitions}
\newtheorem{proposition}[definition]{Proposition}
\newtheorem{example}[definition]{Example}
\newtheorem{comments}[definition]{Some comments}
\newtheorem{corollary}[definition]{Corollary}
\def\square{\Box}
\newtheorem{observation}[definition]{Observation}
\newtheorem{defsobs}[definition]{Definitions and Observations}
\newenvironment{prf}[1]{ \trivlist
\item[\hskip \labelsep{\it
#1.\hspace*{.3em}}]}{~\hspace{\fill}~$\square$\endtrivlist}
\newenvironment{proof}{\begin{prf}{Proof}}{\end{prf}}

\title{A local--global problem for  linear \\ differential equations}
\author{Marius van der Put $\ $ and $\ $ Marc Reversat \\
\footnotesize Department of Mathematics, University of Groningen,
 P.O.Box 800,\\
\footnotesize 9700 AV Groningen, 
The Netherlands, mvdput@math.rug.nl,$\;\;\;$ and\\
\footnotesize Institut de Math\'ematiques de Toulouse, UMR 5219,  Universit\'e Paul 
Sabatier \\
\footnotesize 31602 Toulouse cedex 9, France, email: 
marc.reversat@math.ups-tlse.fr }

\maketitle

{\bf Abstract}.  \begin{footnotesize} An inhomogeneous linear differential equation $Ly=f$ over a global differential field can have a formal solution for each place without having a global solution. The vector space $lgl(L)$ measures this phenomenon. This space is interpreted in terms of cohomology of linear algebraic groups and is computed for abelian differential equations and for regular singular equations.  An analogue of Artin reciprocity for abelian differential equations is given. Malgrange's work on irregularity is reproved in terms 
cohomology of linear algebraic groups.\end{footnotesize}

\section{Introduction} 
The topics:  Elliptic curves $E$ over a number field $K$;
 Drinfeld modules over a field $K$ like $\mathbb{F}_q(t)$; linear differential equations over a differential field $K$, e.g., a finite extension of the differential field $\mathbb{C}(z)$, have many common features. \\
For every place $v$ of $K$ one considers the completion $K_v$. An example 
of a local--global problem is the following. Consider an elliptic curve $E$ over
$K$ and an integer $n>1$. Suppose that $n\cdot y=f$ with $f\in 
E(K)$ has a solution $y_v$ in every $E(K_v)$. Does there
exists a solution $y\in E(K)$? By `folklore' the answer is positive and 
the analogous problem for Drinfeld modules, where the integer $n$ is replaced by 
a non zero element of $\mathbb{F}_q[z]$, has a negative answer (see [H] for both statements).

 Here we consider a differential operator 
 $L=a_n\partial ^n+\cdots +a_1\partial +a_0$, where $\partial$ is the derivation
 on $K$ extending $\frac{d}{dz}$ on $\mathbb{C}(z)$ and $a_n,\dots ,a_0\in K$,
 acting upon $K$. Suppose that the equation $L(y)=f$ with $f\in K$ has a solution
 $y_v$ in every completion $K_v$. Then, in general, there is no solution  in $K$.
 {\it One defines the $\mathbb{C}$-vector space $lgl(L)$} as the kernel of the obvious 
 $\mathbb{C}$-linear map $K/L(K)\rightarrow \prod _v K_v/L(K_v)$. This vector space
 measures this global--local problem. 
 The theme of this paper is the interpretation and the computation of $lgl(L)$.  The main results are: \\
* A formula expressing $lgl(L)$ in terms of the cohomology for differential
 Galois groups acting on the solution space of $L$ and a proof of $\dim lgl(L)<\infty$,\\
* Computation of $H^*(G,V)$ for certain affine group schemes $G$,\\   
* Classification of abelian differential equations and  Artin reciprocity,\\ 
 * Explicit computations of $lgl(L)$ for abelian differential operators and
 for regular singular operators $L$,\\
* A new proof of B.~Malgrange's results on irregularity.\\

The last item requires a precise knowledge of the universal differential Galois group
for the differential field $F=\mathbb{C}(\{z\})$ of the convergent Laurent series.
The multisummation theory of J.-P.~Ramis et al. provides this knowledge.

\section{$K/L(K)$ and $M/\partial M$ as cohomology groups} 
 
$K$ denotes a differential field and let $\partial(f)$ or $f'$ denote the derivative of
$f\in K$. We suppose that its field of constants
$C:=\{f |\ f'=0\}$ is algebraically closed, different from $K$ and has characteristic 0. 
A linear differential equation over $K$ can be written in an operator form
\[(a_n\partial ^n+\cdots +a_1\partial +a_0)y=f\mbox{ with all }a_i,f\in K.\] 
An equivalent formulation is given by a differential module $M=(M,\partial )$
where $M$ is a finite dimensional vector space over $K$ and the additive
operator $\partial :M\rightarrow M$ satisfies $\partial (fm)=f'm+f\partial m$.
(The meaning of the symbols $\partial$ and $'$ will be clear from the context).

\medskip
We recall (see [vdP-S] for details), that for every linear differential equation $M$ (or module) over $K$ there is a differential ring $PVR(M/K)$, called the Picard--Vessiot ring of $M$ over $K$, such that `all solutions' of $\partial m=0$ live in this differential  ring and this ring
has only trivial differential ideals.  The (differential) Galois group of the module $M$ is the linear algebraic group over $C$ consisting of all $K$-linear automorphisms of $PVR(M/K)$, commuting with the differentiation on $PVR(M/K)$. The direct limit of all $PVR(M/K)$
is the universal differential extension $U_K$ of $K$. Its Galois group 
$G_{K,\partial}$ is the affine group scheme, which is the projective limit of the
Galois groups of all differential modules $M$ over $K$.

For a differential operator $L\in K[\partial ],\ L\neq 0$, the solution space of $L$ is the
$C$-vector space  $V(L):=\ker (L,U_K)$ is the solution space of $L$. The action of $G_{K,\partial}$ on $U_K$
leaves $V(L)$ invariant and the restriction of the action to $V$ is the (differential) Galois group of $L$. Similarly, for a differential module $M$ over $K$, the $C$-vector space
$V(M):=\ker (\partial ,U_K\otimes _KM)$ is the solution space of $M$ and the restriction
to $V(M)$ of the natural action of $G_{K,\partial}$ on $U_K\otimes _KM$ is the (differential)
Galois group of $M$. 

\begin{proposition} $H^0(G_{K,\partial},U_K)=K$ and 
$H^i(G_{K,\partial},U_K)=0$ for all $i\geq 1$.
\end{proposition}
\begin{proof} Using that $U_K$ is a direct limit of Picard--Vessiot rings, one 
concludes that it suffices to consider a differential module $M$ over $K$
with Picard--Vessiot ring $L\supset K$ and Galois group $G$. In this case
one has to prove $H^0(G,L)=L^G=K$ and $H^i(G,L)=0$ for all $i\geq 1$.

 The first statement is well known. The affine variety corresponding to $L$
is known to be a $K$-torsor for the linear algebraic group $G$ over $C$. In other
words, there exists a finite Galois extension $K^+\supset K$ such
that $K^+\otimes _KL\cong K^+\otimes _{C}C[G]$. Further the action
of the Galois group ${\rm Gal}(K^+/K)$ on this object commutes with
the action of $G$. 

We recall that the cohomology groups $H^i(G,V)$, where $V$ is any
$G$-module, are the cohomology groups of the Hochschild complex $C^*(G,V)$.  
Moreover, one has $H^i(G,C[G])=0$ for all $i\geq 1$, since
$C[G]=Ind ^G_{\{1\}}C$ is an injective module (see [Jan] for these statements).
It follows that also $H^i(G,K^+\otimes _C C[G])=0$ for $i\geq 1$. Let $i\geq 1$ and let $\xi$ be an
element in $C^i(G,L)$ with $d^i\xi =0$. Then $\xi =d^{i-1}\eta$ for some 
$\eta \in C^{i-1}(G,K^+\otimes L)$. Then $\eta ^*:=\frac{1}{\# Gal(K^+/K)}\cdot 
\sum _{\sigma \in Gal(K^+/K)} \sigma (\eta )$ belongs to $C^{i-1}(G,L)$ and satisfies
$d^{i-1}\eta ^*=\xi$. \end{proof}

\begin{lemma}
Let $L\in K[\partial ],\ L \neq 0$ be a differential operator and $U_K, G_{K,\partial}$ be as before.  The following sequences are exact.
\[0\rightarrow ker(L,U_K)\rightarrow U_K\stackrel{L}{\rightarrow}U_K\rightarrow 0\mbox{ and } \]
\[0\rightarrow ker(L,K)\rightarrow K\stackrel{L}{\rightarrow}K\rightarrow
H^1(G_{K,\partial},ker(L,U))\rightarrow 0.\]
In particular, $K/L(K)$ is canonically isomorphic to $H^1(G_{K,\partial}, V(L))$ and
moreover, $H^i(G_{K,\partial},V(L))=0$ for $i\geq 2$.
\end{lemma}
\begin{proof} For the exactness of the first sequence we have  to show that $Ly=f$ with 
$f\in U_K$ has a solution $y\in U_K$. For any $f\in U_K$ there exists a
$\tilde{L}\in K[\partial ],\ \tilde{L}\neq 0$ such that $\tilde{L}f=0$. Indeed, $f$ lies in some 
$PVR(M)$ and apply now  [vdP-S], Corollary 1.38. The equation $\tilde{L}(L(y))=0$ has all its solutions in $U_K$. Hence, for a suitable solution $y\in U_K$ of this equation one has $Ly=f$.

Taking in the first exact sequence, invariants under $G_{K,\partial}$, one obtains, by using Proposition 2.1, the exact sequence
\[0\rightarrow ker(L,K)\rightarrow K\stackrel{L}{\rightarrow}K\rightarrow
H^1(G,ker(L,U_K))\rightarrow 0.  \] \end{proof}  
 
For a differential module $M$ over $K$ one has a similar result, namely:\\
{\it There is a canonical isomorphism $M/\partial M\rightarrow H^1(G_{K,\partial}, V(M))$
where again $V(M)=\ker (\partial ,U_K\otimes _KM)$ provided with its structure of
$G_{K,\partial}$-module}.  
 
 A direct comparison between modules and differential operators is given by  the 
theorem of the cyclic element which states that any differential module $M$ has the form
 $K[\partial ]/K[\partial ]L$ for some operator $L$. Let $L^*$ denote the adjoint of $L$.
 Then one can verify that $V(M)$ can be identified with $V(L^*)$ and $M/\partial M$ 
 with $K/L^*(K)$.    
 
\section{On cohomology of linear algebraic groups} 

The base field $k$ is supposed to be algebraically closed and to have characteristic 0.
Let $G$ be a linear algebraic group, or, more generally, an affine group scheme, over $k$.  A $G$-module $V$ is a finite dimensional vector space over $k$, provided with an algebraic action of $G$, i.e., a morphism of affine group schemes $\rho :G\rightarrow {\rm GL}(V)$. The cohomology groups $H^*(G,V)$ are defined as the derived functors of $V\mapsto V^G$. The $k$-vector space $H^i(G,V)$ can, as in the case of ordinary group cohomology, be described as the space of all $i$-cocycles 
$f:G\times \cdots \times G\rightarrow V$, divided out by the subspace of the  trivial $i$-cocycles. The only difference is that the $i$-cocycles $f$ are supposed to be morphisms of algebraic varieties (or more generally, of affine schemes). We will also allow
$G$-modules $V$ of infinite dimension, namely  direct limits of finite dimensional 
$G$-modules. Now we collect here (with some comments) the facts and results that we will need in the sequel and refer to [Jan] for the general theory.

\begin{fact} Let $G$ be a reductive affine group scheme (not necessarily connected), then for every $G$-module $V$  one has $H^i(G,V)=0$ for all $i\geq 1$.
\end{fact}
 Indeed, since the characteristic of $k$ is 0, the functor $V\mapsto V^G$ is exact.  
  
\begin{fact} Let $N$ be a closed, normal subgroup of the affine group scheme $G$.
Suppose that $G/N$ is reductive. Let $V$ be a $G$-module. Then $H^*(G,V)$
is canonically isomorphic to $H^*(N,V)^{G/N}$.  
\end{fact}
Indeed, the functor $V\mapsto V^G$ factors as $H^0(G/N,V^N)$ and the functor
$H^0(G/N,-)$ is exact and maps injective objects to acyclic objects. The special case of 3.2, where $N$ is the unipotent radical $R_u(G)$, reduces the computations to the case of (connected) unipotent groups.  

\begin{fact}[the five terms exact sequence]
Let $N$ be a closed normal subgroup of the affine group $G$ and let $V$ be
a $G$-module. Then the exact sequence of five terms reads
\begin{small}
\[0\rightarrow H^1(G/N,V^N)\rightarrow H^1(G,V)\rightarrow H^1(N,V)^{G/N}\rightarrow
H^2(G/N,V^N)\rightarrow H^2(G,V) . \]  
 \end{small}
\end{fact}
This exact sequence is derived  from the spectral sequence $H^a(G/N,H^b(N,V))$ converging to $H^{a+b}(G,V)$.

\begin{remark}  Explicit action of $G/N$ on $H^1(N,V)$.\\ {\rm
Let $V$ be a $G$-module and $N$ a normal closed subgroup of $G$. The action of
 $G/N$ on $H^1(N,V)$ can be made explicit on the level of 1-cocycles. Let 
 $f:N\rightarrow V$ be a 1-cocycle and $g\in G$ then $(gf)$ is the 1-cocycle defined by  
 $(gf)(n)=gf(g^{-1}ng)$. The trivial 1-cocycle  $f(n)=nv-v$ (for a fixed $v\in V$) is mapped to the trivial 1-cocycle $n\mapsto n(gv)-(gv)$. Thus $G$ acts on $H^1(N,V)$. For 
 $m\in N$ and a 1-cocycle $f:N\rightarrow V$ one observes that $(mf)-f$ is the trivial 1-cocycle $n\mapsto nf(m)-f(m)$.  Thus $N$ acts trivially on $H^1(N,V)$.\\ }\end{remark}

\begin{fact} Let $\frak{g}$ be the Lie algebra of a connected affine group scheme
$G$. Any $G$-module $V$ has the structure of a $\frak{g}$-module. The cohomology groups $H^*(G,V)$ are canonically isomorphic to the cohomology
groups $H^*(\frak{g},V)$. 
\end{fact}
{\it Sketch of the proof}. Using  Fact 3.2, we may suppose that $G$ is a connected unipotent
linear algebraic group over $\mathbb{C}$. In particular, $G$ is simply connected. 
Therefore the category of the finite dimensional representations of $G$ and those of its Lie algebra $\frak{g}$ are equivalent. This equivalence extends to an equivalence between the representations of $G$
which are direct limits of finite dimensional representations and those of $\frak{g}$. The latter categories contain enough injective objects and the cohomology groups can be obtained from injective resolutions. \hfill $\square$

We note that  the group $H^i(\frak{g},V)$ can also be described  by $i$-cocycles modulo trivial
$i$-cocycles (see [Jac]). In some cases the computations of the 
$H^*(\frak{g},V)$ are easier than those for $H^*(G,V)$.

\begin{lemma}  Let $t$ be a generator of the Lie algebra of $\mathbb{G}_a$. The
$\mathbb{G}_a$-module $V$, given by $\rho :\mathbb{G}_a\rightarrow {\rm GL}(V)$,
induces a nilpotent map $\rho(t)\in {\rm End}(V)$. Then 
$H^0(\mathbb{G}_a,V)=\ker (\rho(t),V)=\ker (exp(\rho (t))-1,V)$, \\
$H^1(\mathbb{G}_a,V)={\rm coker}(\rho(t),V)={\rm coker}(exp(\rho (t))-1,V)$\\
 and $H^i(\mathbb{G}_a,V)=0$ for $i\geq 2$. \end{lemma}
 \begin{proof} By 3.5, it suffices to compute the cohomology of $V$ as module over
 the Lie algebra of $\mathbb{G}_a$. An obvious calculation of these cohomology groups in terms of cocycles gives the required answer.  The same computation can
 be done in terms of cocycles for the group $\mathbb{G}_a$. \end{proof}

\begin{corollary} Let the group scheme $G$ be topologically generated by
one element $A$. Suppose that the unipotent factor $A_u$ of the Jordan decomposition $A=A_{ss}\cdot A_u$ is non trivial. Then $G=G_1\times \mathbb{G}_a$, where
the first factor is topologically generated by $A_{ss}$ and the second by $A_u$.

Let the $G$-module $V$ be given by $\rho :G\rightarrow {\rm GL}(V)$. Then
 $H^i(G,V)=0$ for $i\geq 2$ and  $H^1(G,V)={\rm coker}(\rho (A_u)-1,V^{G_1})
 ={\rm coker}(\rho (A)-1,V)$.
\end{corollary}
\begin{proof} $H^*(G,V)=H^*(\mathbb{G}_a,V)^{G_1}$ and this equals 
$H^*(\mathbb{G}_a,V^{G_1})$ since $G_1$ commutes with $\mathbb{G}_a$.
Apply now 3.6. The final statement follows from the decomposition of $V$ into
eigenspaces for $\rho(A_{ss})$. \end{proof}

\begin{leeg}  Group schemes with free unipotent generators. \\ {\rm
Let $S$ be some non empty set. One considers tuples $(V,\alpha )$ where
$V$ is a finite dimensional $k$-vector space and $\alpha :S\rightarrow 
{\rm GL}(V)$ maps every $s\in S$ to a unipotent  $\alpha (s)$. This defines an abelian
  category $Unipotent(S)$ which is in an obvious way a Tannakian category. The fibre functor is given by $(V,\alpha )\mapsto V$. Thus $Unipotent(S)$ is isomorphic to the category of the finite dimensional representations of a certain affine group scheme $N$ over $k$. 

Consider an object $(V,\alpha)$ and the Tannakian subcategory $\{\{(V,\alpha )\}\}$ generated by it. The affine group scheme corresponding to this subcategory can
be seen to be the smallest algebraic subgroup $H(V,\alpha )$ of ${\rm GL}(V)$ containing all $\alpha (s)$. The representation $\rho :N\rightarrow {\rm GL}(V)$, corresponding to $(V,\alpha )$, has the property $\rho (N)=H(V,\alpha )$. In particular, if 
$\rho(N)$ is a finite group then $\rho (N)=\{1\}$. Thus $N$ is a connected affine group scheme. Moreover, $N$ is the projective limit of the $H(V,\alpha) $, taken over all objects
$(V,\alpha )$. For a fixed $s\in S$, each $H(V,\alpha )$ contains an element $\alpha (s)$.
The projective limit of the elements $\alpha (s)\in H(V,\alpha )$ can be considered as an
element, again called $s$, in $N$. Thus $\rho (s)=\alpha (s)$ for every object
$(V,\alpha )$ and corresponding representation $\rho :N\rightarrow {\rm GL}(V)$. 
Therefore, we will call  $N$ the {\it group with free unipotent generators $S$}.  Indeed, the elements of
 $S$, seen as elements of $N$, are topological generators, they are unipotent and they have no relations.

 We want to show that for any $N$-module $V$, the $H^*(N,V)$ are the cohomology groups of the complex $0\rightarrow V\rightarrow V^S\rightarrow 0$, where the non trivial
 map is given by $v\mapsto (\rho (s)v-v)_{s\in S}$. A direct proof is maybe possible, however we will prove this using the (pro) Lie algebra of $N$.}\end{leeg}

\begin{leeg} Lie algebras with free nilpotent generators.\\  {\rm
Let $S$ be again a non empty set. The free Lie algebra $F$ over $k$ with generators
$S\subset F$ has the universal property that the representations  $\tau :F\rightarrow {\rm End}(V)$ on vector spaces $V$ over $k$  are in bijection with the maps $S\rightarrow {\rm End}(V)$. One considers now representation $\tau $ such that $V$ is finite dimensional and
every $\tau (s)$ is nilpotent. The image $F(\tau )\subset {\rm End}(V)$ is an algebraic
Lie algebra, because it is generated by nilpotent maps. The corresponding (connected)
algebraic subgroup of ${\rm GL}(V)$ is the smallest algebraic subgroup  containing
all $\exp(\tau (s))$. Define $\rho :S\rightarrow {\rm GL}(V)$ by $\rho (s)=\exp (\tau (s))$.
Then $(V,\rho )$ is an object of $Unipotent(S)$, the Tannakian group of this object
is the smallest algebraic group containing all $\exp(\tau (s))$ and its Lie algebra is
$F(\tau)$.

The projective limit 
$\underset{\leftarrow}{\lim \ }F/\ker{\rho}$, taken over all these $\rho$ will be called the
{\it (pro)-Lie algebra with free nilpotent generators $S$}. It is clear from the above that this pro Lie algebra is the Lie algebra of the affine group scheme $N$ of 3.8. The bijection between the representations of $N$ and those of its (pro)-Lie algebra can be interpreted 
 as  $N$ being simply connected. }\end{leeg}

\begin{proposition} Let $N$ be the affine group scheme with free unipotent generators $S$ and let $V$ be a $N$-module. The $H^*(N,V)$ are the cohomology groups of the complex $0\rightarrow V\rightarrow V^S\rightarrow 0$, where the non trivial
 map is given by $v\mapsto (\rho (s)v-v)_{s\in S}$. 
 \end{proposition}
\begin{proof} Let $Lie(N)$ denote the pro-Lie algebra of $N$. According to 3.5 and the
constructions in 3.8, 3.9, the proposition is equivalent to the statement that the
$H^*(Lie(N),V)$ are the cohomology groups of the complex 
$0\rightarrow V\rightarrow V^S\rightarrow 0$, where the non trivial map is given by
$v\mapsto ( \log (\rho (s))v )_{s\in S}$.\\

For $H^0(Lie(N),V)$, this is obvious. For $H^1(Lie(N),V)$ we use the explicit definition
of a 1-cocycle $f$ (see [Jac]) and conclude that the elements $f(s)\in V$ are arbitrary 
and that they determine $f$ completely. The trivial 1-cocycles are of the form
$f(z)=\tau (z)v$ for a fixed $v\in V$ and where $\tau :Lie(N)\rightarrow {\rm End}(V)$
is induced by $\rho :N\rightarrow {\rm GL}(V)$. This proves the statement for
$H^1(Lie(N),V)$.  The verification of $H^i(Lie(N),V)=0$ for $i\geq 2$ is easy.
\end{proof}

\begin{leeg} Finiteness conditions. \\{\rm
Let the set $S$ be a disjoint union of (non empty) subsets $S_i$ with $i\in I$ (and $I$
an infinite set). Let $Unipotent(\{S_i\}_{i\in I})$ denote the category for which the objects 
are the pairs $(V,\alpha )$ with $V$ a finite dimensional $k$-vector space and 
$\alpha :S\rightarrow {\rm GL}(V)$ such that $\alpha (s)$ is unipotent for every $s\in S$
and there is a finite subset $J$ of $I$ such that $\alpha (s)=1$ for $s\not \in
\cup _{i\in J}S_i$. As in 3.8 and 3.9, this defines an affine group scheme $N$ and
a pro-Lie algebra $Lie(N)$.  The analogue of Proposition 3.10 is:\\

\noindent {\it
Let $V$ be an $N$-module. The $H^*(N,V)$ are the cohomology groups of the complex
$0\rightarrow V\rightarrow V^{(S)}\rightarrow 0$, where $V^{(S)}$ denotes the 
$k$-vector space of the maps $f:S\rightarrow V$ such that there exists a finite subset
$J$ of $I$ with $f(s)=0$ for $s\not \in \cup _{i\in J}S_i$. The non trivial map is again
$v\mapsto (\rho (s)v-v)_{s\in S}$}. \\

The above situation occurs in connection with the Stokes phenomenon (see [vdP-S]), where $S=\{\Delta _{q,d}\ | \ q\in \mathcal{Q},\  d\mbox{ is singular  for } q\}$ is the set of alien derivations. The set $S$ is the disjoint union, over $q\in \mathcal{Q}$, of the
sets $S_q:=\{d\in \mathbb{R} |\ d\mbox{ is singular for }q\}$.  The corresponding 
affine group scheme $N$ is the kernel of the surjective morphism of affine schemes
$G_{analytic}\rightarrow G_{formal}$. We will return to this in Section 5. We note that the finiteness condition  stated in [vdP-S] is slightly wrong. }\end{leeg}

\section{Formal differential equations}
Let $\widehat{F}$ be the differential field $\mathbb{C}((z))$ with derivation $\delta :=z\frac{d}{dz}$. A differential equation or module over $\widehat{F}$ will be called {\it formal}. We recall the explicit descriptions of $U_{\widehat{F}}$ and $G_{formal}:=G_{\widehat{F},\partial}$, slightly extending the one given in [vdP-S]. The {\it aim of this section} is to make both $H^*(G_{formal},V(M))$ (where $V(M)$ is the solution space of $M$) and the canonical isomorphism 
$M/\partial M\rightarrow H^1(G_{formal},V(M))$ explicit.\\

\noindent {\it Description of $U_{\widehat{F}}$.}\\
 Write $\mathbb{C}=A\oplus \mathbb{Q}$ with $A$ a 
$\mathbb{Q}$-vector space. First one introduces the universal Picard--Vessiot ring $U_{rs}$ for the regular singular differential modules over $\widehat{F}$. This ring is 
$U_{rs}:=\overline{\widehat{F}}[\{e(c)\}_{c\in A},\ell ]$ where $\overline{\widehat{F}}$ is the
algebraic closure of $\widehat{F}$ and where the symbols $e(c)$ and $\ell$ satisfy
only the identities $e(c_1+c_2)=e(c_1)\cdot e(c_2)$. The differentiation is given by
$\delta e(c)=c\cdot e(c)$ and $\delta \ell =1$. Then one introduces the set
$\mathcal{Q}=\bigcup _{m\geq 1}z^{-1/m}\mathbb{C}[z^{-1/m}]$ and symbols $e(q)$ for 
$q\in \mathcal{Q}$ satisfying only the identities $e(q_1+q_2)=e(q_1)\cdot e(q_2)$. Further
$\delta e(q)=q\cdot e(q)$. Now $U_{\widehat{F}}:=\oplus _{q\in \mathcal{Q}}U_{\widehat{F},q}$
where $U_{\widehat{F},q}:=U_{rs}e(q)$. This is the universal Picard--Vessiot ring for 
$\widehat{F}$. Put $\mathcal{Q}^+:=\mathcal{Q}\oplus A$ and write $e(q+c):=e(q)\cdot e(c)$
for $q+c\in \mathcal{Q}^+$. Then $U_{\widehat{F}}=\overline{\widehat{F}}[\{e(q)\}_{q\in \mathcal{Q}^+},\ell ]$. \\

\noindent {\it Description of $G_{formal}$.}\\
For convenience, we will identify the affine group scheme $G_{formal}$ with its set 
of $\mathbb{C}$-valued points which is the group of the differential automorphisms of 
$U_{\widehat{F}}/\widehat{F}$.  A special (and very natural) element in this group is the {\it formal monodromy} 
$\gamma$ defined by:\\
(i) $\gamma$ acts on $\overline{\widehat{F}}$ by $\gamma (z^\lambda )=e^{2\pi i \lambda}z^\lambda$ for all $\lambda \in \mathbb{Q}$,\\
(ii) $\gamma e(c)=e^{2 \pi i c}e(c)$ for all $c\in A$,\\
(iii) $\gamma (e(q))=e(\gamma q)$ for all $q\in \mathcal{Q}$,\\
(iv) $\gamma (\ell )=\ell +2\pi i$.\\  

The {\it exponential torus} $T$ (in the terminology of J.-P.~Ramis) is the group 
${\rm Hom}(\mathcal{Q},\mathbb{C}^*)$. An element $h\in T$ acts on $U_{\widehat{F}}$
by $h$ is the identity on $\overline{\widehat{F}}$, on $\ell$ and on the elements $e(c)$.
Further $h(e(q))=h(q)\cdot e(q)$ for all $q\in \mathcal{Q}$. The group $T$ together with the
element $\gamma$ generate topologically $G_{formal}$ (for the Zariski topology). So far we have followed [vdP-S]. The Zariski closure $\overline{<\gamma >}$ of the group 
$<\gamma >$ generated by $\gamma$ is rather big and we prefer to split this group into smaller pieces. For this purpose we introduce more special elements in $G_{formal}$.\\

One decomposes $\gamma$ as a product of commuting automorphisms $\gamma _{ss}\gamma _u$, where $\gamma _{ss}$ has the same definition as $\gamma$ except for (iv) which is replaced by $\gamma _{ss}\ell =\ell$. Further $\gamma _u$ is the identity for the elements
in $\overline{\widehat{F}}$, the elements $e(c)$, $e(q)$ and $\gamma _u\ell =\ell +2\pi i$.
We note that $\gamma _{ss}\gamma _u$ is the Jordan decomposition of $\gamma$ as a product
of a semi-simple element and a unipotent element.

We still want to decompose the semi-simple $\gamma _{ss}$ as a product of commuting elements  
$\gamma _0$ and $\gamma _1$. The direct sum $\mathbb{C}=\mathbb{Q}\oplus A$ yields,
using $c\mapsto e^{2\pi ic}$, a direct product decomposition $\mathbb{C}^*=(\mu _\infty)\times
e^{2\pi iA}$. Now $\gamma _0$ and $\gamma _1$ are the unique semi-simple elements in
$G_{formal}$ with eigenvalues in $\mu_\infty$ and $e^{2\pi iA}$ such that $\gamma _0\gamma _1=\gamma _{ss}$. One verifies that $\gamma _0$ and $\gamma _1$ can also be defined by\\
(i) $\gamma _0 (z^\lambda )=e^{2\pi i\lambda}z^\lambda $, $\gamma _1(z^\lambda )=z^\lambda$
 for all $\lambda \in \mathbb{Q}$;\\
(ii) $\gamma _0(e(q))=e(\gamma _0q)$, $\gamma _1e(q)=e(q)$ for all $q\in \mathcal{Q}$,\\
(iii) $\gamma _0(e(c))=e(c)$, $\gamma _1e(c)=e^{2\pi ic}e(c)$ for all $c\in A$,\\
(iv) $\gamma _0(\ell )=\gamma _1(\ell)=\ell $.\\

 For every element $a\in U_{\widehat{F}}$ there is an integer $n\geq 1$ such that
 $\gamma _0^n(a)=a$. It follows that the algebraic subgroup of $G_{formal}$, generated (topologically) by $\gamma _0$ is $\widehat{\mathbb{Z}}$, the projective limit of the groups 
 $\mathbb{Z}/n\mathbb{Z}$. The algebraic subgroup generated (topologically) by $\gamma _1$
 can be identified with the torus ${\rm Hom}(A,\mathbb{C}^*)$. The algebraic subgroup generated by $\gamma _u$ can be identified with $\mathbb{G}_a$.

Thus $\overline{<\gamma >}$, the Zariski closure of the group generated by $\gamma$, can be identified with 
$\mathbb{G}_a\times {\rm Hom}(A,\mathbb{C}^*)\times \widehat{\mathbb{Z}}$. Moreover, $\gamma$ is the topological generator of $G_{rs}$, i.e., the group of the
differential automorphisms of $U_{rs}/\widehat{F}$, in other words the universal differential  Galois group for the regular singular equations over $\widehat{F}$.

 We extend the exponential torus $T$ to a larger torus
 $T^+={\rm Hom}(\mathcal{Q}^+,\mathbb{C}^*)$. An element $h\in T^+$ acts on 
$U_{\widehat{F}}$ by $h$ is the identity on $\overline{\widehat{F}}$ and $\ell$.
Further $h(e(q)):=h(q)\cdot e(q)$ for $q\in \mathcal{Q}^+$. \\

We conclude that $G_{formal}$ has the form $\mathbb{G}_a\times (T^+\rtimes \widehat{\mathbb{Z}})$. The subgroup $T^+\rtimes \widehat{\mathbb{Z}}$ is reductive and the subgroup
$\mathbb{G}_a$ is the unipotent radical of $G_{formal}$. Further 
$G_{formal}/G^o_{formal}=\widehat{\mathbb{Z}}$ with topological generator $\gamma _0$.\\

\noindent {\it Description of differential modules over $\widehat{F}$.}\\
One associates to a differential module $M$ over $\widehat{F}$ its {\it solution space}:\\ 
$V=V(M):=\ker (\partial ,U_{\widehat{F}}\otimes _{\widehat{F}}M)$ with the following data and relations:\\
(i) a direct sum decomposition $V=\oplus _{q\in \mathcal{Q}}V_q$ where $V_q=\ker (\partial ,U_{\widehat{F},q}\otimes _FM)$,\\
(ii) an element $\gamma _V\in {\rm GL}(V)$ which is the restriction of $\gamma$ on 
$U_{\widehat{F}}\otimes _{\widehat{F}}M$ to the subspace $V$.\\
(iii) $\gamma _V(V_q)=V_{\gamma q}$ for all $q\in \mathcal{Q}$.\\

The functor $M\mapsto V(M):=(V,\{V_q\},\gamma _V)$ defines an equivalence of categories. We note that $\gamma _V$ has a decomposition as a product of commuting elements 
$\gamma _{V,0}\gamma _{V,1}\gamma _{V,u}$ induced by $\gamma =\gamma _0\gamma _1\gamma _u$. \\

\noindent {\it Computation of $H^1(G_{formal},V(M))$}.\\
 By Section 3, $H^i(G_{formal}, V)=H^i(\mathbb{G}_a,W)$ with 
 $W:=V^{T^+\rtimes \widehat{\mathbb{Z}}}$. From the above one sees that $W=V_0^{\gamma _{ss}}$.
 The subspace of $U_{\widehat{F}}$, consisting of the elements invariant under 
$T^+\rtimes \widehat{\mathbb{Z}}$, is $\widehat{F}[\ell ]$. Hence $W=\ker (\partial ,\widehat{F}[\ell ]\otimes _{\widehat{F}}M)$. The action of $\gamma _u$ on $W$ is induced by its action on
$\widehat{F}[\ell ]$, given by $\ell \mapsto \ell +2\pi i$. The group $\mathbb{G}_a$, generated for the
Zariski topology by $\gamma _u$, acts on $W$ by $t \mapsto \gamma  _{u,V}^t$. Hence 
$H^i(G_{formal},V(M))$ can be identified with $\ker (\log \gamma _{u,V},W)$ for $i=0$ and
with ${\rm coker}(\log \gamma _{u,V},W)$ for $i=1$. In particular, the two cohomology groups have the same dimension. This proves Corollary 4.1.\\
 
 \noindent {\it  The explicit canonical isomorphism $can: M/\partial M\rightarrow 
 H^1(G_{formal},V(M))$.}\\
 Any module $M$ is, after a finite extension of $\widehat{F}$,  a direct sum of isotypical 
 submodules modules $M(q)$ with $q\in \mathcal{Q}$. This is just a translation of the decomposition
 $V=\oplus V_q$. Now $\partial$ is bijective on $M(q)$ for $q\neq 0$. Thus for $q\neq 0$, the summand $M(q)$ gives no contribution for the above map $can$. The direct summand $M(0)$, which is the regular
 singular part of $M$, can be represented by a matrix differential operator $z\frac{d}{dz}+B$ with $B$ a constant matrix such that the eigenvalues $\lambda$ of $B$ satisfy $0\leq Re(\lambda )<1$. Only the generalized eigenspace for the eigenvalue $0$ of $B$ can give a contribution to  $M(0)/\partial M(0)$.
This generalized eigenspace is the submodule of $M(0)\approx \widehat{F}^d$ corresponding to $W=V_0^{\gamma _{ss}}$. After restricting to this submodule, $B$ is nilpotent. The kernel and cokernel of $z\frac{d}{dz}+B$ on $\widehat{F}^d$ coincide with the kernel and cokernel of $B$ on
$\mathbb{C}^d$. The fundamental matrix for the equation is
$e^{-B\log z}$. The columns of this matrix form a basis for $\ker (\partial ,\widehat{F}[\ell ]\otimes
_{\widehat{F}}M)$. The differential Galois group is generated by the action of $\gamma _u$ which
multiplies the fundamental matrix by $e^{-2\pi iB}$. Thus $\log \gamma _u$ on $W$ can be identified with $\frac{-1}{2\pi i}B$ on $\mathbb{C}^d$.\\

\begin{corollary}[B.~Malgrange]{\rm [M]}. For any differential module $M$ over $\widehat{F}$ one has
$\dim _\mathbb{C} \ker (\partial ,M) =\dim _\mathbb{C} M/\partial M$.
\end{corollary}

\section{Analytic differential equations}
 $F$ denotes the differential field $\mathbb{C}(\{z\})$, consisting of the convergent
 Laurent series, with derivation $\delta :=z\frac{d}{dz}$. A differential equation or module
 over this field will be called {\it analytic}. We recall and extend results of [vdP-S].\\

 \noindent {\it Description of $U_F$ and $G_{analytic}:=G_{F,\partial}$.}\\
 As in Section 4, $\widehat{F}\supset F$ denotes the
 field of the formal Laurent series. Let $\mathcal{D}\subset \overline{\widehat{F}}$
 denote the $\overline{F}$-subalgebra consising of the elements $h$ that satisfy a linear
 differential equation with coefficients in $\overline{F}$.  Then 
 \[\overline{F}[\{e(q)\}_{q\in \mathcal{Q}^+},\ell ]\subset U_F=\mathcal{D}[\{e(q)\}_{q\in \mathcal{Q}^+},\ell ]\subset U_{\widehat{F}}=\overline{\widehat{F}}[\{e(q)\}_{q\in \mathcal{Q}^+}.\ell ] .\]
{\it The precise structure of the differential algebra $\mathcal{D}$ is unknown}. However the multisummations $\{mult _d\}_{d\in \mathbb{R}}$ in the directions $d$ lead to `locally unipotent' elements of $G_{analytic}$, the Stokes maps (multipliers) for every direction 
$d$. One considers the `locally nilpotent' logarithms $\Delta_d$ of the Stokes maps. These are 
elements of the pro-Lie algebra of $G_{analytic}$. The element $\Delta_d$ is
an $\overline{F}$-linear derivation of $U_F$, commuting with $\delta$, and trivial
on $\overline{F}[\{e(q)\}_{q\in \mathcal{Q}^+},\ell ]$. Thus $\Delta_d$ is determined by its restriction $\Delta_d:\mathcal{D}\rightarrow U_F=\oplus _{q\in \mathcal{Q}}U_{F,q}$. The resulting maps $\Delta_{q,d}:\mathcal{D}\rightarrow U_{F,q}$ are important ingredients for the structure of 
$G_{analytic}$. Put $\{x_i\}_{i\in I}=\{\Delta _{q,d}\ |\  d \mbox{ singular for } q  \}$ and let 
$Lie (N)$ and $N$ denote the pro-Lie algebra and the affine group scheme corresponding to
$\{x_i\}$ defined in Section 3. Then $G_{analytic}=N\rtimes G_{formal}$. We recall that $G_{formal}$ is topologically generated by $T={\rm Hom}(\mathcal{Q},\mathbb{C}^*)$ and the formal monodromy $\gamma$.  The action, by conjugation, of $G_{formal}$ on $N$ induces an action on $Lie(N)$. One can verify that $\gamma (\Delta _{q,d})=\Delta _{\gamma (q),d-2\pi }$ and that, for any $h\in T$, one has 
$h(\Delta _{q,d})=h(q)\Delta _{q,d}$. 

 The group $G_{analytic}/(N\rtimes \mathbb{G}_a)=T^+\rtimes \widehat{\mathbb{Z}}$  is reductive and therefore $N\rtimes \mathbb{G}_a$ is the unipotent radical of 
$G_{analytic}$. The affine group scheme $N\rtimes \mathbb{G}_a$ is topologically generated by  $\{\exp( \Delta _{q,d}) | \ d\mbox{ singular direction for } q\}\cup \{
\gamma _u\}$. The relations between these generators are unknown. The same holds
for their logarithms $\{\Delta _{q,d}\}\cup \{\log \gamma _u\}$.\\
For the group $N\rtimes \overline{< \gamma >}$ with topological generators $\{\exp(\Delta_{q,d})\}\cup \{\gamma \}$ one observes that the only relations are  
$\gamma  \exp (\Delta_{q,d})\gamma ^{-1}=\exp(\Delta _{q,d-2\pi})$. A good translation
in terms of generators of Lie algebras does not seem to exists.\\

\noindent {\it Description of the analytic differential modules}.\\ 
 As in [vdP-S] one can describe a differential module $M$ over $F$ by a structure on
 its space of solutions $V:=V(M)$, namely $V(M)=(V,\{V_q\},\gamma ,\{St_d\})$.
 The maps $St_d$ can be replaced by $st_d=\log St_d$ and $st_d$ can be replaced by its
 components $st_{q,d}$ (with $d$ is a singular direction for $q$). These elements in 
 ${\rm End}(V)$ map each $V_{q'}$ to $V_{q+q'}$ etc. Thus $M$ is represented by a tuple
 $(V,\{V_q\},\gamma ,\{st_{q,d}\})$.\\

 \noindent {\it Computation of $H^*(G_{analytic},V(M))$}.
 The computation uses the formula $H^*(G_{analytic},V(M))=H^*(N\rtimes \mathbb{G}_a,V(M))^{T^+\rtimes \widehat{\mathbb{Z}}}$. First an example.
 
  \begin{example} $M=Fe$ where $\partial (e)=-pe$ and $p=a_kz^{-k}+\cdots +a_1z^{-1}$ with $k\geq 1$ and $a_k\neq 0$. Then $M/\partial M\approx  H^1(G_{analytic},V(M))$ has dimension 
 $k=irr(M)$ (see below for the definition of the {\em irregularity} of $M$). {\rm
 \begin{proof} Let $\rho$ denote the action of $G_{analytic}$ on 
 $V=V(M)=V_p=\mathbb{C}e(p)e$. Then $\rho$ is trivial for $N$ and $\gamma$.
 Further, for $h\in T={\rm Hom}(\mathcal{Q},\mathbb{C}^*)$ the map $\rho(h)$ is
 multiplication by $h(p)$.

 Now $H^1(G_{analytic},V)=H^1(N\rtimes \overline{<\gamma _u >},
 V)^{T\rtimes \overline{<\gamma _{ss}>}}$. The cohomology group
 $H^1(N\rtimes \overline{<\gamma _u>},V)$ identifies with the complex vector space of
the algebraic homomorphism $f:N\rtimes \overline{<\gamma _u >}\rightarrow V$,
such that there are only finitely many $q\in \mathcal{Q}$ for which there is a $d$ with
$f(\exp (\Delta _{q,d}))\neq 0$.

Let $f$ be such  a map and suppose that $f$ is invariant under $T\rtimes \overline{<\gamma _{ss}>)}$.  For $g\in T\rtimes \overline{<\gamma _{ss}>)}$ and 
$m\in N\rtimes \overline{<\gamma _u>}$ one has $f(m)=(gf)(m)$. Since the latter
is $\rho(g^{-1})f(gmg^{-1})$ one has $\rho(g)f(m)=f(gmg^{-1})$. This applied with $g\in T$ and $m=\gamma _u$ or $m=\exp(\Delta _{q,d})$ with $q\neq p$ yields $f(\gamma _u)=0$ and $f(\exp (\Delta _{q,d}))=0$. Now we apply this with  
$m=\gamma _u\exp (\Delta _{p,d})\gamma _u^{-1}$ and
$g=\gamma _{ss}$. Then $f(m)=f(\exp(\Delta_{p,d-2\pi }))$. Using that $f(\gamma _u)= 0$ and $\rho (\gamma _u)=1$ one has $f(m)=f(\exp(\Delta _{p,d}))$.

We conclude that the invariant algebraic homomorphism are described by:
$f(\gamma _u)=f(\exp(\Delta _{q,d}))=0$ for $q\neq p$ and 
$f(exp(\Delta _{p,d}))=f(\exp(\Delta _{p,d-2\pi }))$ for all $d$.
The dimension of the $H^1$ under consideration is therefore equal to the number of the singular directions $d$ modulo $2\pi$ of $p$. This number is easily seen to be $k$. We note that this result coincides with the explicit calculations in [vdP-S], Section 7.3. \end{proof} }\end{example} 
 
 B.~Malgrange has introduced the {\it irregularity of a differential module $M$ over $F$} as follows. Put $\widehat{M}:=\widehat{F}\otimes _F M$. The action of the operator 
 $\partial$ on the exact sequence  $0\rightarrow M\rightarrow \widehat{M}\rightarrow 
 \widehat{M}/M\rightarrow 0$ induces the long exact sequence
 
 \[0\rightarrow \ker (\partial,M)\rightarrow \ker (\partial ,\widehat{M})\rightarrow
 \ker (\partial ,\widehat{M}/M)\]
 \[ \rightarrow M/\partial M\rightarrow \widehat{M}/\partial
 \widehat{M}\rightarrow {\rm coker}(\partial ,\widehat{M}/M)\rightarrow 0 .\]
 
  His results are: ${\rm coker}(\partial ,\widehat{M}/M)=0$, each term in this sequence
  has finite dimension and 
  $\dim \ker(\partial ,\widehat{M})=\dim \widehat{M}/\partial \widehat{M}$ (see
  Corollary 4.1).  The {\it irregularity of $M$}, $irr(M)$, is defined as the dimension of $\ker (\partial ,\widehat{M}/M)$. Using cohomology we will reprove 
 Malgrange's results and compute  $irr(M)$.

 \begin{proposition}[B.~Malgrange]{\rm [M]}.  Let $M$ be a differential module over $F$ with solution space $V(M)=(V,\{V_q\},\gamma ,\{st_{q,d}\})$. Then $\partial$ is surjective on 
 $\widehat{M}/M$ and $irr(M)=\sum _q\deg _{z^{-1}}q\cdot \dim V_q$. 
 \end{proposition}
 \begin{proof}  The map $M/\partial M\rightarrow
  \widehat{M}/\partial \widehat{M}$ identifies with the restriction map $R:$
  $H^1(N\rtimes G_{formal},V)\rightarrow H^1(G_{formal},V)$, where $V$ is the solution space of $M$.  The five term exact sequence (with $G=G_{analytic}=N\rtimes G_{formal}$) reads
 \begin{footnotesize}
\[0\rightarrow H^1(G_{formal},V^N)\rightarrow H^1(G_{analytic},V)\rightarrow 
H^1(N,V)^{G_{formal}}\rightarrow
H^2(G_{formal},V^N)\rightarrow \cdot \cdot \]  
\end{footnotesize}
The term   $H^2(G_{formal},V^N)$ is zero. The description of $H^1(N,V)$ by 1-cocycles
modulo trivial 1-cocycles and the description of the pro Lie algebra $Lie(N)$ yield the exact sequence
$0 \rightarrow V/V^N\rightarrow V^{(\{\Delta _{q,d}\} )}\rightarrow H^1(N,V)\rightarrow 0,$  Taking invariants for the action of $G_{formal}$ one finds the exact sequence
\[0\rightarrow H^0(G_{formal},V/V^N)\rightarrow H^0(G_{formal},V^{(\{\Delta _{q,d}\} )})\rightarrow \]\[
H^0(G_{formal},H^1(N,V))\stackrel{\alpha _1}{\rightarrow} H^1(G_{formal},V/V^N)\rightarrow 
H^1(G_{formal},V^{(\{\Delta _{q,d}\} )}). \]

We will compute the terms of this exact sequence.\\

\noindent 
(1) {\it Claim}:  $Irr:= H^0(G_{formal},V^{(\{\Delta _{q,d}\} )})$ has dimension 
 $\sum \deg_{z^{-1}}q\cdot \dim V_q$.
Further $H^1(G_{formal},V^{(\{\Delta _{q,d}\} )})=0$ and so $\alpha _1$ is surjective.\\

\noindent 
 $V^{(\{\Delta _{q,d}\} )}$ consists of the maps $f:\{\Delta _{q,d}\}\rightarrow V$ (with the property that only 
 finitely many $q$'s have a $d$ with $f(\Delta _{q,d})\neq 0$). The action of $\gamma$
 on $V^{(\{\Delta _{q,d}\} )}$ is given by $(\gamma f)(\Delta _{q,d})=\rho(\gamma ^{-1})f(\Delta _{\gamma(q),d-2\pi}))$. The map $f\mapsto (\gamma f)-f$ is seen to be surjective. It follows that 
 $H^1(\overline{<\gamma >},V^{(\{\Delta _{q,d}\} )})=0$ and that \\
 $H^0(\overline{<\gamma >},V^{(\{\Delta _{q,d}\} )})$ consists of those
 maps  satisfying $f(\Delta _{q,d})=f(\Delta _{\gamma(q),d-2\pi })$. Taking now the invariants under $T$ one finds that $H^1(G_{formal},V^{(\{\Delta _{q,d}\} )})=0$ and 
 $H^0(G_{formal},V^{(\{\Delta _{q,d}\} )})$ consists of the maps satisfying $f(\Delta _{q,d})\in V_q$ and  $f(\Delta _{q,d})= f(\Delta _{\gamma(q),d-2\pi})$.

 Write $V=V_0\oplus V_{q_1}\oplus \cdots \oplus V_{q_r}$ where $0,q_1,\dots ,q_r$
 are distinct elements of $\mathcal{Q}$ and the $V_{q_i}\neq 0$. We allow the possibility
 $V_0=0$.  From the above it follows that $Irr:=H^0(G_{formal},V^{(\{\Delta _{q,d}\} )})$ has dimension $\sum _i \deg_{z^{-1}}q_i\cdot \dim V_{q_i}$.

\bigskip

\noindent (2) {\it Claim}: The map $R: H^1(N\rtimes G_{formal},V)\rightarrow H^1(G_{formal},V)$ is surjective.\\

First we consider the long exact sequence
\[0\rightarrow H^0(G_{formal},V^N)\rightarrow H^0(G_{formal},V)\rightarrow 
H^0(G_{formal},V/V^N)\]
\[\stackrel{\alpha _2}{\rightarrow}H^1(G_{formal},V^N)\rightarrow 
H^1(G_{formal},V)\rightarrow H^1(G_{formal},V/V^N)\rightarrow 0 .\]
We observe that $H^0(G_{formal},V^N)=\ker (\partial ,M)$ and that
$H^0(G_{formal},V)$ is equal to $\ker (\partial ,\widehat{M})$.

The morphism from the exact sequence
\[0\rightarrow H^1(G_{formal},V^N)\rightarrow H^1(G_{analytic},V)\rightarrow
H^0(G_{formal},H^1(N,V))\rightarrow 0\]
to the exact sequence
\[ 0\rightarrow H^1(G_{formal},V^N)/im(\alpha _2) \rightarrow H^1(G_{formal},V)\rightarrow H^1(G_{formal},V/V^N)\rightarrow 0,\]
is defined by the three maps $\alpha _3, \ R, \ \alpha _1$, where $\alpha _3$
is induced by the identity map $H^1(G_{formal},V^N)\rightarrow H^1(G_{formal},V^N)$.
 The surjectivity of $\alpha _1$ implies that $R$ is surjective.

\bigskip

\noindent (3) By (2) one has an exact sequence 
$0\rightarrow im(\alpha _2)\rightarrow \ker (R)\rightarrow \ker (\alpha _1)\rightarrow 0$.
Combining this with the results, observations and identifications of (1) and (2) one finds the Malgrange's exact sequence
\[0\rightarrow \ker(\partial ,M)\rightarrow \ker(\partial ,\widehat{M})\rightarrow Irr
\rightarrow M/\partial M\rightarrow \widehat{M}/\partial \widehat{M}\rightarrow 0\ .\]
  \end{proof}

\section{Global differential equations}

Now $K$ is a finite extension of the differential field $\mathbb{C}(z)$. The universal $U_K$
and $G_{K,\partial}$ are far from known. For every place $v$ one has a canonical embedding $K\subset K_v$. 
We note that $K_v$ is isomorphic to  $\mathbb{C}((t))$, where $t$ is a local parameter  
for $v$. The differentiation has the form $a\cdot \frac{d}{dt}$ for some non zero element $a$
(which will not be of importance for our problems). 

The embedding $K\subset K_v$ can be extended to 
 an injection $U_K\rightarrow U_{K_v}$ of differential rings. This injection is not unique, however its image is. Thus the above arrow is unique up to a differential automorphism
 of $U_K/K$, i.e., an element of $G_{K,\partial}$. For every place $v$ we make a choice for
 $U_K\rightarrow U_{K_v}$.  This arrow induces an injective morphism $G_{K_v,\partial}\rightarrow G_{K,\partial}$ (which is unique up to conjugation by an element of $G_{K,\partial}$).
 
Let $L\in K[\partial ],\ L\neq 0$ and $f\in K$. The equation $L(y)=f$ has at a place $v$
a solution $y_v\in K_v$ if $L$ is regular at $v$ and $f$ has no pole at $v$. Thus for all but
finitely many places there exists a solution $y_v$. Hence there is a well defined map
of $\mathbb{C}$-vector spaces $K/L(K)\rightarrow \oplus _vK_v/L(K_v)$. The kernel
of this map will be denoted by $lgl(L)$. This vector space measures the failure of the
the local--global principle for differential equations.

 For a differential module $M$ over $K$ one defines in a similar way the $\mathbb{C}$-vector space 
$lgl(M):=\ker M/\partial M\rightarrow \oplus _v(K_v\otimes M)/\partial (K_v\otimes M)$.

Using the above and Section 2 one concludes that
\[lgl(M)=\ker (H^1(G_{K,\partial},V(M))\rightarrow \oplus _vH^1(G_{K_v,\partial },V(M)))\]
and similarly for a differential operator $L$.

\begin{theorem} The $\mathbb{C}$-vector spaces $lgl(L)$ and $lgl(M)$ have finite dimension. 
\end{theorem}
\begin{proof} The statements for $L$ and $M$ are equivalent. We start with a differential
module $M$ over $K$. As $K$ is a finite extension of $F:=\mathbb{C}(z)$ we may 
view $M$ as a differential module over $F$. The term $M/\partial M$ does not change.
Consider a place $w$ of $F$ and the places $v_1,\dots ,v_r$ of $K$ above $w$. Then
$F_w\otimes _FM$ can be identified with $\oplus _{i=1}^rK_{v_i}\otimes _KM$. This implies
that $lgl(M)$ does not change if one considers $M$ as differential module over $F$.
Thus we may suppose that $K=\mathbb{C}(z)$.

 Now we consider $L=a_n\partial ^n+\cdots +a_1\partial +a_0$ with all $a_i\in K=\mathbb{C}[z]$
 and $g.c.d.(a_n,\dots ,a_0)=1$. Consider an equation $L(y)=f$ with $f\in K$
 which has a solution $y_v\in K_v$ for every place $v$. For each  place $v\neq \infty$
 we write this solution as $[y_v]+r_v$, where $[y_v]=\sum _{i\geq 1}\frac{c_i}{(z-v)^i}$ is the principle part of $y_v$.
 Only finitely many $[y_v]$ are non zero and $g:=f-\sum _{v\neq \infty} L([y_v])$ lies
 in $\mathbb{C}[z]$. Therefore it suffices to consider equations $L(y)=g$ with
 $g\in \mathbb{C}[z]$. 
 
  One observes that the operator 
 $L:\mathbb{C}[z]\rightarrow \mathbb{C}[z]$ has a finite dimensional cokernel. This implies that $lgl(L)$ has finite dimension. More precisely, $lgl(L)$ is the kernel of the
 obvious map $\mathbb{C}[z]/L(\mathbb{C}[z])\rightarrow \oplus _{{\rm all}\ v}
 \mathbb{C}(z)_v/L(\mathbb{C}(z)_v)$.  \end{proof}

 \section{Abelian differential equations}
 As in Section 2, $K$ is a differential field with an algebraically closed field of constants 
 $C\neq K$ of characteristic 0. A differential module $M$ (or an operator in $K[\partial ]$) is called {\it abelian} if the differential Galois group $Gal(M)$ is abelian. The corresponding Picard--Vessiot extension $L\supset K$ is also called {\it abelian} and its (differential) Galois group is denoted by $Gal(L/K)$.

 Any abelian linear algebraic group $G$ is a product of copies of $\mathbb{G}_a,\ \mathbb{G}_m$
and finite cyclic groups. The `{\it additive part}' $G^+$ of $G$ is (by definition) a product of copies of $\mathbb{G}_a$ and the `{\it multiplicative part}' $G^{\star}$ of $G$ is (by definition) a product of copies of $\mathbb{G}_m$ and finite cyclic groups. \\

\noindent 
 A linear algebraic group $G$ of additive type is described as $G=Spec(C[W])$
 where $W$ is a finite dimensional $C$-vector space, $C[W]:=\oplus _{n\geq 0}sym ^nW$ and where the comultiplication $m$ is given by $m(w)=(w\otimes 1)+(1\otimes w)$ for all $w\in W$. A group scheme of additive type $G$ is (by definition) the projective
 limit of linear algebraic groups of additive type. It follows that it has also a
 description as $G=Spec(C[W])$ but now with $W$ any $C$-vector space. The set $G(C)$ of the $C$-valued points of $G$ is clearly identified with ${\rm Hom}_C(W,C)$.
 The set ${\rm Hom}(G,\mathbb{G}_a)$ of the morphisms can be  identified with $W$   
 as follows. Fix a presentation $\mathbb{G}_a=C[t]$ with $m(t)=(t\otimes 1)+(1\otimes t)$. Giving a morphism $G\rightarrow \mathbb{G}_a$ is equivalent to giving a $w\in W$, image of $t$ for the map $C[t]\rightarrow C[W]$. By abuse of language, we will write
${\rm Hom}_C(W,C)$ for the group scheme $Spec(C[W])$ of additive type.\\

A linear algebraic group $G$ of multiplicative type can be described as 
$G=Spec(C[A])$ with $A$ is a finitely generated abelian group, $C[A]$ the group ring
of $A$ over $C$ and with comultiplication $m(a)=a\otimes a$ for all $a\in A$. A group
scheme of multiplicative type is again defined as $G=Spec(C[A])$, but now $A$ is an
arbitrary abelian group. Now $G(C)={\rm Hom}_\mathbb{Z}(A,C^*)$ and 
${\rm Hom}(G,\mathbb{G}_m)=A$. By abuse of language, we will write
${\rm Hom}_\mathbb{Z}(A,C^*)$ for the group scheme $Spec(C[A])$ of multiplicative
type. Note that ${\rm Hom}(G_1,G_2)=0$
 if $G_1$ is of additive type and $G_2$ is of multiplicative type or visa versa.\\
 
For any differential field $L$ one defines the $C$-vector subspace \\
$d(L):=\{f'\ |\ f\in L\}\subset L$
 and the subgroup $dLog(L):=\{\frac{f'}{f}\ |\ f\in L^*\}\subset L$.

\begin{theorem} 
Let $L$ be an abelian Picard-Vessiot extension of $K$, then we have natural maps, induced by the inclusion
$K\subseteq L$,
	\[\alpha_{L/K} :\frac{K}{dLog(K)}\rightarrow \frac{L}{dLog(L)} \ \ {\text{ and }} \ \ 
	\beta_{L/K} :\frac{K}{d(K)}\rightarrow \frac{L}{d(L)}.\]
\begin{enumerate}
	\item $\ker (\alpha_{L/K}) $ is a ${\mathbb Z}$-module of finite type and ${\rm Hom}_{\mathbb Z}(\ker (\alpha_{L/K}) ,C^{\ast})$ is considered as a linear algebraic group (of multiplicative type).
	\item $\ker (\beta_{L/K}) $ is a finite dimensional $C$-vector space. ${\rm Hom}_{C}(\ker (\beta_{L/K}) , C)$ is considered as a linear algebraic group (of additive type).
	\item There are canonical isomorphism of linear algebraic groups
	\[\psi_{L/K}^{\star}: Gal(L/K)^{\star}\rightarrow {\rm Hom}_{\mathbb Z}(\ker (\alpha_{L/K}) ,C^{\ast})  \mbox{ and}\]
\[	\psi_{L/K}^{+}: Gal(L/K)^+\rightarrow {\rm Hom}_{C}(\ker (\beta_{L/K}) , C). \]

\end{enumerate}
\end{theorem}

The {\it definitions of $\psi_{L/K}^{\star}, \psi _{L/K}^+$} are the following. Take 
$\sigma \in Gal(L/K)$, $\bar{f}\in \ker (\alpha_{L/K})$ and $\bar{g}\in \ker (\beta_{L/K})$, images
of $f\in K$ and $g\in K$. There are elements $F\in L^*,\ G\in L$ such that $\frac{F'}{F}=f$
and $G'=g$. Then $\sigma (F)=cF$ for some $c\in C^*$ and $\sigma (G)=G+d$ for some
$d\in C$. Then
 \[\psi _{L/K}^{\star}(\sigma)(\bar{f})=c\in C^*\mbox{ and }\psi _{L/K}^+(\sigma )(\bar{g})=d\in C. \] 
It is easily verified that the above definitions do not depend on the choices for $f,F,g,G$ and
that the maps $\psi_{L/K}^\star,\ \psi _{L/K}^+$ are homomorphisms of algebraic groups.
`Canonical' means that for abelian Picard--Vessiot extensions $K\subset L_1\subset L_2$ one 
has the obvious rules for composing the various $\psi ^\star ,\ \psi ^+$. 

\begin{proof}
 The direct elementary proof that $\psi ^\star,\ \psi ^+$ are isomorphisms, is somewhat long. We prefer to use a result of  [vdP-R], namely  Lemma 1.1, which reduces the general case to the 
 cases $Gal(L/K)\in \{\mathbb{G}_m, C_n,\mathbb{G}_a\}$, where $C_n$ denotes the cyclic
 group of order $n>1$. For these cases we provide elementary proofs.\\
 (i) $Gal(L/K)=\mathbb{G}_m$. Then $L$ is the Picard--Vessiot extension of an equation
 $y'=fy$ with $f\in K$, such that $y'=mfy$ has, for every integer $m\geq 1$, in $K$ only the trivial solution $y=0$. Then $L$ is the transcendental extension $K(t)$ with $t'=ft$ and $c\in Gal(L/K)=C^*$ maps $t$ to $c t$.\\
The kernel of $\psi _{L/K}^\star : K/dLog(K)\rightarrow L/dLog(L)$ consist of the elements 
$h+dLog(K)$ such that there exists $y\in L^*$ with $\frac{y'}{y}=h$. The obvious elements $y\in L^*$ with $\frac{y'}{y}\in K$ are $at^n$ with $a\in K^*$ and $n\in \mathbb{Z}$. They produce
$\mathbb{Z}\bar{f}$ in the kernel of $\psi _{L/K}^\star$. We want to show that there are no other
elements $y\in L^*$ with $\frac{y'}{y}\in K$. We may write $y=\frac{p(t)}{q(t)}$
with $p(t),q(t)\in K[t]$ relatively prime, $p(t),q(t)$ both monic and not divisible by $t$. Then 
$\frac{y'}{y}=\frac{p(t)'}{p(t)}-\frac{q(t)'}{q(t)}\in K$ implies that $\frac{p(t)'}{p(t)}$ and $\frac{q(t)'}{q(t)}$ both belong to $K$. Now $p(t)=t^n+\cdots +p_1t+p_0,\  p_0\neq 0$. Then 
$p(t)'=p_{n-1}'t^{n-1}+\cdots +p_1't+p_0'+ f(nt^n+\cdots +p_1t)=a\cdot p(t)$ for some $a\in K$.
This identity generates equalities $a=nf,\ p_0'=ap_0$ and $p_0'=nfp_0$ contradicts the assumption on the equation $y'=fy$ if $n>0$. Hence $p(t)=1$ and similarly $q(t)=1$.
 Thus we find that the kernel of $\psi _{L/K}^\star$ is
$\mathbb{Z}\bar{f}$. The map $Gal(L/K)=C^*\rightarrow {\rm Hom}(\mathbb{Z}\bar{f},C^*)$
is obviously an isomorphism. \\
 
 \noindent (ii) $Gal(L/K)=C_n$. Then $L$ is the Picard--Vessiot extension of an equation $y'=fy$
with $f\in K$ and $n$ is minimal such that $y'=nfy$ has a non zero solution $y_0\in K$.
Then $L=K(t)$ with $t^n=y_0$ and $Gal(L/K)$ acts by multiplying $t$ by $n$th roots of unity.
The proof that the kernel of $\psi _{L/K}^\star$ is the cyclic group $\mathbb{Z}\bar{f}$ of order $n$ is similar to the one of case (i).\\

\noindent (iii) $Gal(L/K)=\mathbb{G}_a$. Then $L=K(t)\neq K$ where $t$ satisfies a differential equation of the form $t'=g$ with $g\in K$. An element $d\in Gal(L/K)=\mathbb{G}_a=C$ maps $t$ to $t+d$.  An element $h+d(K)$ lies in the kernel of $\psi _{L/K}^+$ if and only if there exists
an element $H\in L$ with $H'=h$. The kernel clearly contains the $C$-subspace generated by
$g+d(K)$.  It suffices to show that the kernel contains no more elements. Write 
$H=\frac{p(t)}{q(t)}+r(t)$ with relatively prime $p(t),q(t)\in K[t]$, $q(t)$ monic, 
$\deg p(t)<\deg q(t)$ and $r(t)\in K[t]$. Suppose that  $H'=r(t)'+\frac{p(t)'q(t)-p(t)q(t)'}{q(t)^2}=h
\in K$. Then 
  \[p(t)'q(t)-p(t)q(t)'=(h-r(t)')q(t)^2.\]
 If $q(t)\neq 1$, then, by comparing  the degrees, one finds $r(t)'=h$. The same holds if $q(t)=1$. We may write  $r(t)=r_dt^d+\cdots +r_1t$, since
the the constant term of $r(t)$ is of no importance. Further $r_d\neq 0$. 
Then $r(t)'=r_d't^d+\cdots +r_1't+ g(r_ddt^{d-1}+\cdots +r_1)=h$. This implies
$r_d'=0$ and thus $r_d=c\in C^*$. For $d>1$ one finds the contradiction  $r_{d-1}'+gdc=0$.
 Hence $r(t)=c t$ and $h=cg$. This finishes the computation.  \end{proof}

Theorem 7.1 admits the following corollary, which can be interpreted as an Artin correspondence for abelian Picard--Vessiot extensions.

\begin{corollary}  
Let $K^{\mathrm{ab}}_{\mathrm{diff}}$ be the maximal abelian  Picard--Vessiot extension of $K$, then its Galois group $Gal (K^{\mathrm{ab}}_{\mathrm{diff}}/K)$ satisfies the following isomorphisms of affine group schemes over $C$
	\[\psi_K ^\star:Gal (K^{\mathrm{ab}}_{\mathrm{diff}}/K)^\star
	\rightarrow {\rm Hom}_{\mathbb Z}(\frac{K}{dLog(K)} ,C^{\ast})\mbox{ and }\]
	\[\psi_K ^+:Gal (K^{\mathrm{ab}}_{\mathrm{diff}}/K)^+
	\rightarrow {\rm Hom}_{C}(\frac{K}{d(K)} , C).\]
	Write $\psi _K: Gal (K^{\mathrm{ab}}_{\mathrm{diff}}/K)\rightarrow
	{\rm Hom}_{\mathbb Z}(\frac{K}{dLog(K)} ,C^{\ast})\times {\rm Hom}_{C}(\frac{K}{d(K)} , C)$
	for $\psi _K^\star \times \psi _K^+$.\\
     
     \noindent $\psi _K$ induces a bijective correspondence between the  Picard--Vessiot 
    extensions  of finite type $K\subset L\subset  K^{\mathrm{ab}}_{\mathrm{diff}}$   
	and the  pairs $(Z, V)\subset (K/dLog(K), K/d(K))$, with $Z$ a subgroup of finite type,
	$V$ a $C$-vector subspace of finite dimension.
	\end{corollary}

	We note that $\psi _K^\star$ is the projective limit of the $\psi _{L/K}^\star$ taken over
	all Picard--Vessiot extensions of finite type $L\supset K$ contained in 
	$K^{\mathrm{ab}}_{\mathrm{diff}}$ (and  similarly for $\psi _K^+$).  
	The correspondence can be described as follows.\\
	One associates to $(Z,V)$ the Picard--Vessiot field $L=(K^{\mathrm{ab}}_{\mathrm{diff}})^H$,
	where  $H$ is the kernel of the restriction map $Gal (K^{\mathrm{ab}}_{\mathrm{diff}}/K)\rightarrow  {\rm Hom}_{\mathbb{Z}}(Z,C^*)\times {\rm Hom}_C(V,C)$. \\
This can be made even more explicit by giving $Z$ generators 
$\bar{a}_1,\dots ,\bar{a}_r$ over $\mathbb{Z}$ and $V$ generators $\bar{b}_1, \dots ,\bar{b}_s$ over $C$. Then $L$ is the Picard--Vessiot extension for the set of equations
\[y_i'=a_iy_i\mbox{ for }i=1,\dots ,r \mbox{ and } z_j'=b_j\mbox{ for }j=1,\dots ,s.\]
 Conversely, $(\ker (K/dLog(K)\rightarrow L/dLog(L)), \ker (K/d(K)\rightarrow L/d(L)))$ is 
the pair	associated to $L$.\\
	
We remark that the explicit presentation of the universal Picard--Vessiot ring and its automorphism
group for the category of the abelian differential modules over $K$, as given in [vdP-R], also provides a proof for 7.1 and 7.2.  \\

 \section{Computation of $lgl(L)$ for abelian $L$} 
 
 We consider an abelian differential module $M$ over a finite extension $K$ of
 $\mathbb{C}(z)$ and compute
 the space $lgl(M)$. It suffices to do this for an indecomposable $M$ and this reduces the general case to the cases:\\
 (1) $M$ of dimension 1 and trivial.\\
 (2) $M$ has dimension 1 and differential Galois group $\mathbb{G}_m$.\\
 (3) $M$ has dimension 1 and differential Galois group $C_m$ with $m>1$.\\
 (4) $M$ of dimension $n$ has differential Galois group $\mathbb{G}_a$.\\ 

\subsection{ The case of a trivial $L$}

We consider a differential field $K$ which is a finite extension of $\mathbb{C}(z)$ 
and $L=\frac{d}{dz}$. The solution space $V=\mathbb{C}1$ and $lgl(L)$ is the kernel 
of the map $H^1(G_{K,\partial},V)\rightarrow \oplus _vH^1(G_{K_v,\partial},V)$. Since these groups $G_*$ act trivially on $V$, the term $H^1$ coincides with the morphism
$G_* \rightarrow \mathbb{G}_a$. A morphism factorizes over the additive factor of the 
abelianized group $(G_*)_{ab}^+$.  For the group $G_{K,\partial }$ the additive factor
(see Corollary 7.2) can be written as ${\rm Hom}_C(\Omega _K /dK,C)$, where $\Omega _K$ is the differential module for $K/\mathbb{C}$. Further 
${\rm Hom}((G_{K,\partial})_{ab}^+,\mathbb{G}_a)$ identifies with $\Omega _K/dK$.
Similarly for the groups $G_{K_v,\partial }$.

 Let $X$ be the curve associated to $K$. The points $x$ of $X$ correspond to the 
 $v$'s. Further $\Omega _K$ identifies with the space of all meromorphic differential forms on $X$.  We recall the exact sequence, studied in [vdP-R], Section 1.3
 \[0\rightarrow H^1_{DR}(X,\mathbb{C})\rightarrow \Omega _K/dK\rightarrow
 \oplus _{x\in X}\Omega _x/d(K_x)\rightarrow \mathbb{C}\rightarrow 0\ .\]
This implies the following result.

\begin{proposition} $lgl(\frac{d}{dz})$ for the finite field extension 
$K\supset \mathbb{C}(z)$ associated to a curve $X$ 
over $\mathbb{C}$ of genus $g$, is canonically isomorphic to the $2g$-dimensional 
vector space $H^1_{DR}(X,\mathbb{C})$.
\end{proposition}

\subsection{The $\mathbb{G}_m$ case}

 $K$ is again a finite extension of $\mathbb{C}(z)$, corresponding to a curve $X$ over 
 $\mathbb{C}$. Let the point $p\in X$ have local analytic parameter $t$. 
  One writes, as before, $K_p$ for the completion of $K$ at the valuation induced by $p$.
  Thus $K_p=\mathbb{C}((t))$. The subfield of the convergent Laurent series 
  $\mathbb{C}(\{t\})$ is denoted by $K_p^{an}$. Further 
  $\Omega _{K_p}:=\mathbb{C}((t)) dt$ and $\Omega _{K_p^{an}}:=\mathbb{C}(\{t\})dt$.
  
 The operator $L$ has the form $\frac{d}{dz}-f$ with $f\in K^*$ such that
  $dy =m\cdot yfdz$ has, for  any integer $m\geq 1$, in $K$ only the trivial solution 
  $y=0$. One considers the map, again called $L$, $y\in K\mapsto dy-yfdz\in \Omega _K$. Let $H\subset \Omega _K$ be the subspace
  consisting of the elements $\omega$ such that there exists for every point $p$ a formal local  solution $y\in K_p$ of $L(y)=\omega$. Then $lgl(L)$ is the cokernel of $L:K\rightarrow H$.

  A point $p$ is {\it regular} for $L$ if there exists a $g\in K_p, \neq 0$ with $L(g)=0$. In fact this $g$ lies in $K^{an}_p$.   Further $p$ is regular if and only if $fdz$ has at most a pole of order 1 at $p$ and $Res_p(fdz)\in \mathbb{Z}$. For a regular point $p$ the map $L:K_p\rightarrow \Omega_{K_p}$, the kernel and cokernel have dimension 1. The same holds for  $L:K^{an}_p\rightarrow \Omega _{K^{an}_p}$. This is seen be writing $y=gh$ with $h\in K_p$. Then $L(y)=L(gh)=gdh$.The kernel of $L$ is $\mathbb{C}g$ and $\omega$ lies in the image of $L$ if and only if $Res_p(g^{-1}\omega )=0$. 
  
 $p$ is a {\it regular singular} if $fdz$ has a pole of order 1 at $p$ and 
 $Res_p(fdz)\not \in \mathbb{Z}$. In this case the map $L:K_p\rightarrow \Omega _{K_p}$ is bijective and the same holds for $L:K_p^{an}\rightarrow \Omega _{K^{an}_p}$.
 
 A point $p$ is {\it irregular singular} if $fdz$ has a pole of order $d+1$ with $d>0$. The integer $d$ is the {\it irregularity} $irr _p$ of $p$. In this case the map $L:K_p\rightarrow \Omega _{K_p}$ is bijective. Further $L:K_p^{an}\rightarrow \Omega _{K^{an}_p}$ is injective and its cokernel $Irr(p,L)$ has dimension $d=irr_p$.\\
 
 Let $Sol(L)$ denote the subsheaf of the sheaf of meromorphic functions $M$ on $X$,
defined by $Sol(L)(U)=\{y\in M(U)\ |\ L(y)=0\}$ for any open $U\subset X$. It is a sheaf of            $\mathbb{C}$-vector spaces. By assumption $H^0(X,Sol(L))=0$.

\begin{lemma} $\dim H^1(X,Sol(L))=2g-2+\#S$ and $H^2(X,Sol(L))=0$.
\end{lemma} 
\begin{proof}
Let $S\subset X$ denote the set of singular points of $L$. {\it First we suppose 
$S\neq \emptyset$}.  
The restriction of the sheaf $\mathcal{L}=Sol(L)$ to $X^*:=X\setminus S$ is locally isomorphic to the constant sheaf $\underline{\mathbb{C}}$. It is given by a non trivial homomorphism of the fundamental group $\pi _1:=\pi _1(X^*)\rightarrow {\rm GL}(V)$ where $V$ is a 1-dimensional vector space. Then $H^i(X^*,\mathcal{L})$  equals the cohomology group $H^i(\pi _1 ,V)$ for all $i$. The group $\pi _1$ is free on $r:=2g-1+\#S$ generators. Let $t_1,\dots ,t_r$ denote these free generators.
The action of $t_i$ on $V$ is multiplication by some $\alpha _i\in \mathbb{C}^*$. The cohomology groups $H^*(\pi _1,V)$ are the cohomology groups of the complex $0\rightarrow V\rightarrow V^r\rightarrow 0$, where the non trivial arrow is defined by $v\mapsto (t_iv-v)_{i=1,\dots ,r}$. Since some $\alpha _i\neq 1$
one has $H^i(\pi _1,V)=0$ for $i=0,2$ and $\dim H^1(\pi _1,V)=r-1$. \\

{\it We claim  that $H^1(X,\mathcal{L})\rightarrow H^1(X^*,\mathcal{L})$ is an
isomorphism}. Let $U$ be the disjoint union of small disks around the points of $S$ and let
$U^*=U\setminus S$. The Mayer--Vietoris sequence for the covering $\{X^*,U\}$ yields the exact sequence
\[0\rightarrow H^1(X,\mathcal{L})\rightarrow H^1(X^*,\mathcal{L})\oplus H^1(U,\mathcal{L})
\rightarrow H^1(U^*,\mathcal{L})\rightarrow \cdots . \]
Consider a small disk $D_p$ around a point $p\in S$ and let $D_p^*=D_p\setminus \{p\}$. 
We can identify $D_p$ with $D:=\{z\in \mathbb{C}|\ |z|<1\}$, $D_p^*$ with 
$D^*:=\{z\in \mathbb{C}|\ 0<|z|<1\}$ and the restriction of $\mathcal{L}$ with the kernel of the
morphism $M:O_D\rightarrow O_D$ (here $O_D$ denotes the sheaf of the holomorphic 
functions on $D$), given by $y\mapsto zy'-ay$ with any $a\in \mathbb{C}\setminus \mathbb{Z}$. The sequence $0\rightarrow \ker M\rightarrow O_D\rightarrow O_D\rightarrow 0$ is exact.
One verifies that $M:H^0(E,O_D)\rightarrow H^0(E,O_D)$ is bijective for $E=D$ and $E=D^*$.
Further $H^i(E,O_D)=0$ for $E=D$ and $E=D^*$ and $i=1,2$. This implies that 
$H^i(E,\ker M)=0$ for $E=D$ and $E=D^*$ and $i=1,2$. In this way we have verified that  
$H^1(U,\mathcal{L})=H^1(U^*,\mathcal{L})=0$ and this proves the claim. Moreover $H^2(X,\mathcal{L})=0$.

\bigskip
 
Now we consider {\it the case $S=\emptyset$}. The sheaf $\mathcal{L}=Sol(L)$
is locally isomorphic to the constant sheaf $\underline{\mathbb{C}}$. However,
$\mathcal{L}$ is not equal to the constant sheaf since we have supposed that 
the equation $L(y)=0$ has in $K$ only the solution $y=0$. In particular, the genus $g$ of $X$
is $\geq 1$.  Take any point $p\in X$. Then $H^1(X\setminus \{p\},\mathcal{L})$ is isomorphic
to $H^1(\pi ,V)$, where $\pi$ is the fundamental group of $X^*:=X\setminus \{p\}$, free on $2g$
generators and $V$ is the 1-dimensonal vector space such that the restriction of
$\mathcal{L}$ to $X\setminus \{p\}$ corresponds to a non trivial action of $\pi$ on $V$. 
Thus $\dim H^1(\pi ,V)=2g-1$. Let $D$ be a small disk around $p$ and put 
$D^*=D\setminus \{p\}$. The Mayor-Vietoris sequence for the covering $\{X^*,D\}$ of $X$ yields
an exact sequence
\[0\rightarrow H^1(X,\mathcal{L})\rightarrow H^1(X^*,\mathcal{L})\oplus H^1(D,\mathcal{L})\rightarrow H^1(D^*,\mathcal{L})\rightarrow H^2(X,\mathcal{L})\rightarrow 0\ .\]
Now $H^1(D,\mathcal{L})=0$  and $H^1(D^*,\mathcal{L})$ has dimension 1. According to Lemma 8.3, the space  
$H^2(X, \mathcal{L})$ is {\it  dual} to $H^0(X,\mathcal{L}^*)$ and therefore 0. 
It follows that  $\dim H^1(X,\mathcal{L})=2g-2$.  \end{proof}

For the following, well known, result we could not find a reference.

\begin{lemma} Let $\mathcal{M}$  be a sheaf of $\mathbb{C}$-vector spaces on $X$, locally isomorphic to
the constant sheaf $\underline{\mathbb{C}}^k$. Then $\sum _{i=0}^2\dim H^i(X,\mathcal{M}) =k(2-2g)$.
Moreover, there is a canonical isomorphism $H^2(X,\mathcal{M})\rightarrow H^0(X,\mathcal{M}^*)^*$.
\end{lemma}
\begin{proof} The definition of the dual $\mathcal{M}^*$ is rather obvious. Consider the exact sequence 
$0\rightarrow \mathcal{M}\rightarrow \mathcal{M}\otimes _{\underline{\mathbb{C}}}O_X\rightarrow
\mathcal{M}\otimes _{\underline{\mathbb{C}}}\Omega _X\rightarrow 0$ obtained by tensoring
the exact sequence $0\rightarrow \underline{\mathbb{C}}\rightarrow O_X\rightarrow \Omega _X
\rightarrow 0$ with $\mathcal{M}$ over $\underline{\mathbb{C}}$.

This induces an exact sequence of cohomology groups above $X$
\[ 0\rightarrow H^0(\mathcal{M})\cdots H^1(\mathcal{M}\otimes O_X)\stackrel{A}{\rightarrow} H^1(\mathcal{M}\otimes \Omega _X)\rightarrow
H^2(\mathcal{M})\rightarrow 0\ .\]
The vector bundle $\mathcal{M}\otimes O_X$ has rank $k$. Its degree is 0, since the line bundle $\Lambda ^k(\mathcal{M}\otimes O_X)$ admits
a connection without singularities.  This implies the formula. By Serre duality,  $A$  is the dual of 
$ H^0(\mathcal{M}^*\otimes O_X)\stackrel{B}{\rightarrow} H^0(\mathcal{M}^*\otimes \Omega _X)$. Now $(\ker B)^*\cong {\rm coker }A$ yields  the required duality
$H^2(\mathcal{M})^*\cong H^0(\mathcal{M}^*)$. \end{proof}

One defines the skyscraper sheaf $\mathcal{Q}$ on $X$ by the exact sequence of sheaves
\[0\rightarrow Sol(L)\rightarrow M\stackrel{L}{\rightarrow} \Omega _{mer}\rightarrow \mathcal{Q}\rightarrow 0\ ,\]
where $\Omega _{mer}$ denote the sheaf of the meromorphic differential forms on $X$.  
The sheaf $\mathcal{Q}$ is the sheaf $\oplus _{p\in X} E_p$. Each term $E_p$ denotes a 
sheaf with zero stalk at the points $q\neq p$  and its stalk at $p$ is the finite dimensional vector space  $\Omega _{K_p^{an}}/L(K_p^{an})$. By definition, a section of $\mathcal{Q}$ above
an open set $U$ is an element of $\prod _{p\in U}\Omega _{K_p^{an}}/L(K_p^{an})$ which has
a discrete support.

Let the sheaf $\mathcal{H}$ denote the image of $L$. The exact sequences 
\[0\rightarrow Sol(L)\rightarrow M\rightarrow \mathcal{H}\rightarrow 0 \mbox{ and }\] 
\[0\rightarrow \mathcal{H}\rightarrow \Omega _{mer}\rightarrow \mathcal{Q}\rightarrow 0\]
induce long exact sequences for their cohomology on $X$. The first one yields, in combination with Lemma 8.4, $H^i(X,\mathcal{H})=0$ for $i=1,2$ and the exact sequence
$0\rightarrow K\rightarrow H^0(X,\mathcal{H})\rightarrow H^1(X,Sol(L))\rightarrow 0$.

The second one yields the exact sequence 
\[ 0\rightarrow H^0(X,\mathcal{H})\rightarrow 
\Omega _K\rightarrow \oplus _{p\in X}\Omega _{K_p^{an}}/L(K_p^{an})\rightarrow 0 \ .\]

 The space $H_0:=H^0(X,\mathcal{H})$ consists of the differential forms  
 $\omega \in \Omega _K$ such that at every point $p$ there exists an element $y\in K^{an}_p$
 with $L(y)=\omega$. Let  $H\subset \Omega_K$ denote the subspace of the $\omega$
 such that for  every point $p$ there exists a $y\in K_p$ with $L(y)=\omega$. We recall that
 $lgl(L)$ is the cokernel of $L:K\rightarrow H$.

 The last exact sequence implies that {\it $R:\Omega _K\rightarrow \oplus _{p\in Irr}Irr(p,L)$ is surjective}. Here $Irr$ denotes the set of the irregular singular points and we recall that 
 $Irr(p,L)=\Omega _{K_p^{an}}/L(K_p^{an})$.   It follows that the dimension of $lgl(L)$ is
 the sum of $\sum _{p\in Irr}irr_p$ and the dimension of the cokernel of $L:K\rightarrow H_0$.
 The latter can be identified with $H^1(X, Sol(L))$. Thus we have proved the following.
 
\begin{theorem} $\dim lgl(L)=2g-2+\#S+\sum _{p\in Irr}irr_p$
\end{theorem}

\begin{example} The regular singular, order 1, operator on $\mathbb{P}^1$ {\rm
\[L=\frac{d}{dz}+\sum _{j=1}^r\frac{\lambda _j}{z-p_j} \mbox{ with }
\lambda _j\not \in \mathbb{Z} \mbox{ for all }j\ . \]

\noindent $S=\{p_1,\dots ,p_r\}$ if $\sum \lambda _j\in \mathbb{Z}$ and otherwise
$S=\{p_1,\dots ,p_r,\infty \}$. According to Theorem 8.3,  $\dim lgl(L)$ is $r-2$ 
in the first case  and $r-1$ in the second one.\\

\noindent {\it An explicit calculation of $lgl(L)$ using the proof of Theorem {\rm 6.1.}}\\
 $L$ is replaced by the operator
 \[L^+= \prod _{j=1}^r(z-p_j)L= \prod _{j=1}^r(z-p_j)\frac{d}{dz}+q(z).\]
Then $lgl(L)=lgl(L^+)$ and the latter is equal to
the kernel of 
\[\mathbb{C}[z]/L^+(\mathbb{C}[z])\rightarrow \oplus _{p\in \mathbb{P}^1}
\mathbb{C}(z)_p/L^+(\mathbb{C}(z)_p).\]
The map $\mathbb{C}[z]/L^+(\mathbb{C}[z])\rightarrow \mathbb{C}(z)_p/L(\mathbb{C}(z)_p)$
is the zero map except when $p=\infty$ and $\sum \lambda _j\in \mathbb{Z}$. 
In verifying that  $\dim {\rm coker}( L^+,\mathbb{C}[z])=r-1$,
one has to consider separately the cases: $\sum \lambda _j=-n$ with $n>0$ an integer,
$\sum \lambda _j=0$ and $\sum \lambda _j$ 
is not an integer $\leq 0$.  \hfill $\square$ }\end{example}

\begin{example} The irregular singular  operator 
$L=\frac{d}{dz}-z^5$ on $\mathbb{P}^1$ satisfies $\dim lgl(L)=5$. {\rm 
Indeed, $g=0,\ S=\{\infty \},\  irr_\infty =5$. The latter is verified by writing a multiple of $L$
in terms of the parameter $t=z^{-1}$  as $\frac{d}{dt}+t^6$.

The explicit method of Theorem 6.1 and the observation that $L$ is bijective on
$\mathbb{C}(z)_\infty$ yield that $lgl(L)$ is equal to the cokernel of $L$
on $\mathbb{C}[z]$. The latter has as basis (the images of) 
$1,z,z^2,z^3,z^4$. \hfill $\square$ }\end{example}

\begin{example} $L=z^5\frac{d}{dz}+1+az^4$ with $a\not \in \mathbb{Z}$ acting on 
$\mathbb{P}^1$ has $\dim lgl(L)=4$.\\
{\rm  Indeed, $g=0,\ S=\{0,\infty \}, irr _0=4$. Now the method of Theorem 6.1. The cokernel of $L$ on $\mathbb{C}[z]$ has a basis of representatives $\{1,z,z^2,z^3\}$. The equation $Ly=f$ with $f$ a polynomial of degree  $\leq 3$, has a local solution at $\infty$. Hence $\dim lgl(L)=4$. \hfill $\square$ }\end{example}

\begin{example} $K=\mathbb{C}(z)[y]$ with $y^2=(z-a_1)\cdots (z-a_{2g+2})$. The trivial differential module $K$ with $\partial 1=0$ satisfies $\dim lgl(K)=2g$.\\
{\rm This follows at once from Proposition 8.1. Another way to calculate $lgl(K)$ is to view
$K$ as the two-dimensional differential module $\mathbb{C}(z)1\oplus \mathbb{C}(z)y$ with
$\partial 1=0$ and $\partial y= (\sum _{j=1}^{2g+2}\frac{1/2}{z-a_j})y$. The first factor has
$lgl=0$ and for the second factor, the Example 8.4
yields $\dim lgl =2g$.}\hfill $\square$ \end{example}

\subsection{The $C_m$ case with $m>1$}
This case is rather similar to the case $\mathbb{G}_m$. The only new point is that there are
only regular singularities, since the differential Galois group is finite. With the same 
notations as above, the result is therefore

\begin{proposition} $\dim lgl(L)=2g-2+\#S$.
\end{proposition}

 \subsection{ The $\mathbb{G}_a$ case}
  The action of $\mathbb{G}_a$ on the solution space is given
 by $t\mapsto e^tN$, where $N$ is a nilpotent $n\times n$ matrix consisting of one
 Jordan block. Using that $K$ is a $C_1$-field and that the differential
 Galois group is connected, one finds a matrix differential equation for $M$ of the form
 $\frac{d}{dz}-fN$ for some non zero $f\in K$ (see [vdP-S], Corollary 1.32).
 We note that this differential equation  is regular singular. The {\it singular points} $p\in X$, where $X$ denotes the curve associated to $K$,
 are precisely the points where the residue of $fdz$ is not zero. Let $S\subset X$ denote the set of the singular points. If $S$ is not empty, then $S$ contains at least 2 points since the sum of the residues of $fdz$ is zero. If $S=\emptyset$, then $fdz$ could be exact.
 In that case we are in the situation of Subsection 8.1 and therefore
  {\it we assume that $fdz$ is not exact}. \\
 
 We want to compare the cokernel of 
 $\frac{d}{dz}-fN$ acting upon $K^n$ with the cokernels of the same operator acting upon all $K_p^n$ for $p\in X$.   \\
 
 \noindent {\bf Computation for $n=2$}. We consider the operator
 \[L({y_1\choose y_2}):=(\frac{d}{dz}-fN){y_1\choose y_2}={y_1'-fy_2\choose y_2'}:K^2\rightarrow K^2 
,\mbox{ or equivalently} \]
 \[L((y_1, y_2 )=(dy_1-fy_2dz, dy_2):\  K^2\rightarrow  \Omega _K^2\ .\]
 Let $H\subset \Omega_K^2$ be the complex subspace consisting of the elements
 $(\omega_1, \omega_2)$ such that there is a formal local solution at each point $v$ of
 $\mathbb{P}^1$. Then $lgl(\frac{d}{dz}-fN)$ is the cokernel of the map
 $L:K^2\rightarrow H$. 
 
 We recall that $X$ is the curve associated to $K$.
 Consider the $\mathbb{C}$-linear map $R:H\rightarrow H^1_{DR}(X,\mathbb{C})$, 
 given by $(\omega _1, \omega _2)\mapsto \overline{\omega _2}$,  where 
 $\overline{\omega _2}$ is the image of $\omega _2$ in $H^1_{DR}(X,\mathbb{C})$.
 This is well defined since all the residues of $\omega _2$ are 0.\\
 {\it Now we investigate the image of $R$. }\\
 
 {\it Suppose that $S$ is not empty}. Take a $\omega _2$ (with all residues 0) representing
 a given element in $H^1_{DR}(X,\mathbb{C})$. We have to produce an $\omega _1$ such that 
 $(\omega _1,\omega _2)\in H$. We will use the existence of a meromorphic differential
 form $\omega$ for any prescription of its residues $Res_p(\omega )=a_p$ such that
 almost all $a_p$ are 0 and $\sum a_p=0$.

 For any point $p\not \in S$, such that $\omega _1$ and $fdz$ have no pole at $p$,
 we take  $a_p=0$. For a point $p\not \in S$, such that $\omega _1$ has a pole at $p$ or 
 $fdz$ has a pole at $p$, we define $a_p$ by $a_p+Res_p(y_2fdz)=0$ (where locally at $p$ one has
 $dy_2=\omega _1$). For a singular point $p\in S$ we take a local solution $y_2$
 of $dy_2=\omega _2$. Then, as we may change $y_2$ into $y_2+c$ for any $c\in \mathbb{C}$,
 we have that for any choice of $a_p$, there is a constant $c$ such that the residue of 
 $a_p+ Res_p( y_2fdz +cfdz)=0$. Since $S$ is not empty, we can choose  the last values
 of $a_p$ such that $\sum a_p=0$. Any $\omega _1$ with these residues satisfies
 $(\omega _1,\omega _2)\in H$. {\it Thus $R$ is surjective}.\\
 
 Consider the map $L:K^2\rightarrow H_0:=\{(\omega _1,\omega _2)\in H|\ 
 \overline{\omega _2} =0\}$. Dividing $H_0$ by $L(\{(0,y_2)| y_2\in K\})$, we have to compute the cokernel of $L:K\oplus \mathbb{C}\rightarrow 
 H_{00}$ with $H_{00}:=\{(\omega _1,\omega _2)\in H|\ \omega _2=0\}$.   
  Any $(\omega _1,0)\in H_{00}$ is mapped to the image of 
$(Res _v(\omega _1))_{v\in S}$ in the space  
\[\{(a_v)_{v\in S}\in  \mathbb{C}^S|\ \sum a_v=0\}/\mathbb{C}(Res_v(fdz))_{v\in S}.\]
This map is surjective. The kernel of this map consists of the $\omega _1$ for which
there exists a constant $c$ such that $\omega _1+cfdz$
has all residues 0 (we note that $c$ is unique). This leads to the statement that the cokernel of
$L:K\oplus \mathbb{C}\rightarrow H_{00}$ has dimension $ (-2+\#S)+\dim H^1_{DR}(X,\mathbb{C})$ and the cokernel of $L:K^2\rightarrow H$ has dimension $2g+(-2+\#S)+2g$. \\

 {\it Suppose that $S=\emptyset$}. Let $M$ and $\Omega _{mer}$ be the sheaves on $X$ (for the ordinary complex topology) of the meromorphic 
 functions and the meromorphic differential forms. Both sheaves have $H^i=0$ for $i\geq 1$. One considers the morphism 
 $L: M ^2\rightarrow \Omega ^2_{mer}$  defined by $(y_1,y_2)\mapsto (dy_1-y_2fdz,dy_2)$. Let $\mathcal{H}\subset \Omega _{mer}^2$ be the image
 of $L$. Then $H^0(X,\mathcal{H})$ consists of the pairs $(\omega _1,\omega _2)\in \Omega _K^2$ such that the equation
 $(dy_1-y_2fdz,dy_2)=(\omega _1,\omega _2)$  has everywhere a local solution. The sheaf $\mathcal{L}$ is defined by the exact sequence
 \[0\rightarrow \mathcal{L}\rightarrow M^2\rightarrow \mathcal{H} \rightarrow 0 \ .\]
 The sheaf $\mathcal{L}$ is locally isomorphic to the constant sheaf $\underline{\mathbb{C}}^2$. Taking the cohomology above $X$ one finds the exact sequence
 \[ 0\rightarrow  H^0(X,\mathcal{L})\rightarrow K^2\rightarrow H^0(X,\mathcal{H})\rightarrow H^1(\mathcal{L})\rightarrow 0\ .\]
  Thus $H^1(\mathcal{L})$ identifies with $lgl(\frac{d}{dz}-fN)$.
 It is easily seen that $H^0(\mathcal{L})$ and $H^0(\mathcal{L}^*)$ have both dimension 1. Lemma 8.3 implies that 
 $\dim H^1(\mathcal{L}) =4g-2$.

\begin{proposition} Suppose that $fdz$ is not exact and that $n=2$. Then the dimension of
$lgl(\frac{d}{dz}-fN)$ is $4g-2 +\# S$, where $S$ consists of the singular points
of $\frac{d}{dz}-fN$, i.e.,  the points $p$ with $Res_p(fdz)\neq 0$.
\end{proposition}

  \begin{theorem}  Suppose that $fdz$ is not exact. Consider $L=(\frac{d}{dz}-fN)$, where 
  $N\in {\rm End}(\mathbb{C}^n)$ is a nilpotent matrix with one Jordan block. $S$ is the set
  of singular points of $L$, i.e., the points $p$ with $Res_p(fdz)\neq 0$.
  Then $lgl(L)$ has dimension $2g+(n-1)\cdot (2g-2+\#S)$.
  \end{theorem}   
  
 \begin{proof} For $S=\emptyset$ one easily verifies that the above method, explained for $n=2$,  holds for any $n\geq 2$. 
 For $S\neq \emptyset$ the proof uses induction w.r.t. $n$. For notational convenience we only consider $n=3$.

 $lgl(L)$ is the cokernel of the map $L:K^3\rightarrow H\subset \Omega _K^3$, given by 
 \[L(y_1,y_2,y_3)=(dy_1-y_2fdz,dy_2-y_3fdz,dy_3)\ .\]
  As before $H$ is the subspace of 
 $\Omega _K^3$ consisting of the tuples $(\omega _1,\omega _2,\omega _3)$ for which 
 there are formal local solutions at every point $p\in \mathbb{P}^1$. As before we consider
 the map $R: H\rightarrow H^1_{DR}(X,\mathbb{C})$ which sends a tuple to 
 $\overline{\omega _3}$, the image of $\omega _3$ in $H^1_{DR}(X,\mathbb{C})$. We claim
 that $R$ is surjective.
 
   Let $\omega _3$ with $Res_p(\omega _3)=0$ be given. As in the proof of Proposition 8.2,
 there exists an $\omega _2$ such that the pair of equations 
 $dy_2=y_3fdz +\omega _2,\ dy_3=\omega _3$ has a formal local solution at each point
 of $\mathbb{P}^1$. The reason is that the residues of $\omega _2$ at the points of $S$ can
 be arbitrarily chosen, with the only restriction taking that the sum of all residues is 0.
 Now that $\omega _2$ is chosen we want to produce a $\omega _1$ such that the pair
 of equations $dy_1=y_2fdz +\omega _1,\ dy_2=y_3fdz +\omega _2$ has a formal local 
 solution at each point of $\mathbb{P}^1$. By construction, the last equation has a solution
 $y_2$, unique up to a constant. Using this constant one can prescribe the residues of 
 $\omega _1$ at the points of $S$, with the only restriction that the sum of all residues is 0.
 {\it This shows that $R$ is surjective}.\\
 
 \noindent 
 Now $lgl(L)$ is the direct sum of $H^1_{DR}(X,\mathbb{C})$ and the cokernel of
 $L_0:K^3\rightarrow H_0$ where $H_0\subset H$ consists of the tuples where $\omega _3$
 is exact. After dividing by the subspace $\{L_0(0,0,y_3)\ |\ y_3\in K\}$ we have
 to calculate the cokernel of the map $L_{00}:K\oplus K\oplus \mathbb{C}\rightarrow H_{00}$,
 given by $(y_1,y_2,y_3)\mapsto (dy_1-y_2fdz,dy_2-y_3fdz,0)$ and where
 $H_{00}\subset H$ consists of the tuples with $\omega _3=0$. \\
 
 Consider $(\omega _1,\omega _2,0)\in H_{00}$.  The equation 
 $dy_2=y_3fdz+\omega _2$ is formally solvable at each point $p$ for a suitable 
 $y_3\in \mathbb{C}$ (depending on $p$) if and only if $Res_p(\omega _2)=0$ for every 
 $p\not \in S$. As in the proof of Proposition 8.2, this induces a linear map
 $R_{00}:H_{00}\rightarrow \{(a_p)_{p\in S}\in \mathbb{C}^S| \ \sum a_p=0\}$, given by
 $(\omega _1,\omega _2,0)\in H_{00}$ maps to $(Res_p(\omega _2))_{p\in S}$. The image of
 $fdz$ in the space is non zero. We conclude that the dimension of the cokernel of 
 $L_{00}$ is equal to $-2+\#S\ +$ the dimension of the cokernel of 
 $L_{000}:K\oplus K\oplus \{0\}\rightarrow H_{000}$, where $H_{000}\subset H_{00}$ consists
 of the tuples with such that $Res_p(\omega _2)=0$ for all $p$. By the case $n=2$, the dimension
 of the cokernel is $4g-2+\#S$. Thus the cokernel of $L$ has dimension 
 $(2g-2+\#S)+(4g-2+\#S)$.    \end{proof}
 
  \section{Regular singular differential equations}
  Now we consider a differential operator $L$  on $\mathbb{P}^1$ which has only {\it  regular
  singularities}. The aim is to calculate $lgl(L)$.  One can represent $L$ as a global connection
  $\nabla:\mathbb{C}(z)^m\rightarrow \Omega _{\mathbb{C}(z)/\mathbb{C}}^m$. Let $Mer$ denote  the  sheaf of the meromorphic functions and $\Omega _{mer}$ the sheaf of the
  meromorphic differential forms on $\mathbb{P}^1$. One defines the sheaves $\mathcal{L}$
  and $\mathcal{Q}$ on $\mathbb{P}^1$ by the exact sequence of sheaves
  \[0\rightarrow \mathcal{L}\rightarrow Mer ^m\stackrel{\nabla}{\rightarrow}\Omega _{mer}^m\stackrel{R}{\rightarrow} \mathcal{Q}\rightarrow 0\ .\]
  
  An element of $lgl(L)$ is represented by an $\omega \in \Omega _{\mathbb{C}(z)/\mathbb{C}}^m$ having the property that the equation $\nabla (y)=\omega$ has at every point $P$ of $\mathbb{P}^1$  a formal solution. The assumption that $L$ has only regular singular points implies that this formal solution is in fact meromorphic at the given point and therefore lies in the stalk
  $Mer ^m_P$, or equivalently $R(\omega  )=0$.  \\

 {\it First we analyze the sheaf $\mathcal{Q}$}. Let $t$ be a local parameter at $P$
 (i.e., $t=z-c$ or $t=z^{-1}$). Locally at an open disk $D$ around $P$ the map 
 $\nabla$ has the form $\nabla (y)=dy +A_Py\frac{dt}{t}$, where $A_P$ is a constant
 matrix such that the real part of each eigenvalue of $A_P$ lies in $[0,1)$. We note
 that $A_P=0$ if $P$ is not a singular point of $L$. For a singular point $P$ of
 $L$ we write $t_0(P)$ for the dimension of the kernel of $A_P$. We write 
 $Res_P:\Omega _{mer,P}\rightarrow \mathbb{C}^m_P$ for the residue map and
 we write $Res (\omega )=\sum _PRes_P(\omega )[P]$. The latter is a section of the skyscraper sheaf $\oplus _P\mathbb{C}^m_P$ on $\mathbb{P}^1$.

 An element  $y\in Mer ^m_P$ can be written as 
$y=\sum _{n>>-\infty }y_nt^n$ with all  $y_n\in \mathbb{C}^m$ and an element
$\omega \in \Omega _{mer,P}^m$ as $\sum _{n>> -\infty}\omega _nt^ndt$ with
all $\omega _n\in \mathbb{C}^m$. 
 
 Now $\nabla (y)=\sum _{n>>-\infty} (n+A_P)y_nt^{n-1}dt=\sum _{n>>-\infty}\omega _nt^ndt$ has a solution if and only if $\omega _{-1}=Res_P(\omega )$ lies in the image of $A_P$. It follows that $\mathcal{Q}$ is the skyscraper sheaf  
$\oplus _P(\mathbb{C}^m_P/A_P\mathbb{C}^m_P)$ and $\mathcal{Q}_P$ has
dimension $t_0(P)$. Moreover, $t_0(P)$ is the dimension of the space of 
solutions of $\nabla (y)=0$, locally at the point $P$.
 
 \smallskip
 
 The above exact sequence of sheaves is an acyclic resolution of $\mathcal{L}$ since the sheaves $Mer$ and $\Omega _{mer}$ have `trivial' cohomology on every open subset of $\mathbb{P}^1$. Thus the cohomology groups of the complex
 \[0\rightarrow \mathbb{C}(z)^m\stackrel{\nabla}{\rightarrow}\mathbb{C}(z)^mdz
 \stackrel{R}{\rightarrow}\oplus _P\mathcal{Q}_P\rightarrow 0\]  
 can be identified with $H^*(\mathbb{P}^1,\mathcal{L})$. The very definition of
 $lgl(L)$ yields that $lgl(L)\cong H^1(\mathbb{P}^1,\mathcal{L})$.
 
 \smallskip
  {\it We will 
 need a formula for the dimension of $H^2(\mathbb{P}^1,\mathcal{L})$} (in other terms, the dimension of the cokernel of $R$ in the above complex). The image of 
 $Res:\mathbb{C}(z)^mdz\rightarrow \oplus _P\mathbb{C}^m_P$ consists of the
 elemens $v=(v_P)_P$ with $\sum v_P=0$. Let $p_1,\dots ,p_r$ denote the
 singular points of $L$ and $A_1,\dots ,A_r$ the corresponding constant matrices
 as above. Then one can verify that the cokernel of $R$ is isomorphic to 
 $\mathbb{C}^m/(\sum _{i=1}^rA_i(\mathbb{C}^m))$. Here we have identified all
 $\mathbb{C}^m_P$ with a single vector space $\mathbb{C}^m$. These identifications are not explicit and we can only conclude that 
 $\dim H^2(\mathbb{P}^1,\mathcal{L})<m$ if $r>0$  and that 
 $H^2(\mathbb{P}^1,\mathcal{L})=0$ if some 
 $A_i$ is invertible. The latter is equivalent to $t_0(p_i)=0$ and is again equivalent
 to the statement that the equation $\nabla (y)=0$ has no local solution $y\neq 0$
 in a neighborhood of $p_i$. 
 
 \bigskip
 
   {\it Next, we will compute for a singular point $p_i$ and a small disk $X_i$ around  $p_i$, the dimensions of $H^*(X_i,\mathcal{L})$ and 
   $H^*(X_i^*,\mathcal{L})$, where   $X_i^*=X_i\setminus \{p_i\}$}. \\
   Let $O_{hol}$ and $\Omega _{hol}([p_i])$ denote the sheaves on $X_i$ of the
   holomorphic functions and the differential forms having at most a pole of order 1
   at $p_i$. The connection $\nabla :O_{hol}^m\rightarrow \Omega _{hol}([p_i])$
   has, as before, the form $\nabla (y)=dy +A_iy\frac{dt}{t}$.  
    Define the sheaf $\mathcal{P}$ by the exact sequence (above $X_i$)
    \[0\rightarrow \mathcal{L}\rightarrow O_{hol}^m\stackrel{\nabla}{\rightarrow}
    \Omega _{hol}([p_i])^m\rightarrow \mathcal{P}\rightarrow 0\ .\]
   Then $\mathcal{P}$ is a skyscraper sheaf with at most one non zero stalk, namely $\mathcal{P}_{p_i}=\mathbb{C}^m/A_i(\mathbb{C}^m)$. The above sequence is an acyclic resolution of the restriction of $\mathcal{L}$ to $X_i$.
   Taking global sections on $X_i$ and $X_i^*$ one easily finds that the spaces
   $H^0(X_i,\mathcal{L}),\ H^0(X_i^*,\mathcal{L}),H^1(X_i^*,\mathcal{L})$
   have dimension $t_0(p_i)$ and that the other cohomology groups are 0. 
 
 \bigskip

 Now we start {\it a computation of $\chi (\mathbb{P}^1):=\sum (-1)^i\dim H^i(\mathbb{P}^1,\mathcal{L})$}.  
Write $X_0=\mathbb{P}^1\setminus \{p_1,\dots ,p_r\}$ and let $X_i$ denote, as before, a small
disk around $p_i$ for $i=1,\dots ,r$. Further $X_i^*:=X_i\setminus \{p_i\}$ and 
$Y:=\cup _{i=1}^rX_i$ and $Y^*:=\cup _{i=1}^rX_i^*$.   The Mayer-Vietoris exact sequence for the covering $X_0\cup Y$ of  $\mathbb{P}^1$ reads
\[0\rightarrow H^0(\mathbb{P}^1,\mathcal{L})\rightarrow H^0(X_0,\mathcal{L})\oplus 
H^0(Y,\mathcal{L})\rightarrow H^0(Y^*,\mathcal{L}) \]
\[  \rightarrow H^1(\mathbb{P}^1,\mathcal{L})\rightarrow H^1(X_0,\mathcal{L})\oplus 
H^1(Y,\mathcal{L})\rightarrow H^1(Y^*,\mathcal{L}) \]
\[  \rightarrow H^2(\mathbb{P}^1,\mathcal{L})\rightarrow H^2(X_0,\mathcal{L})\oplus 
H^2(Y,\mathcal{L})\rightarrow H^2(Y^*,\mathcal{L})\rightarrow 0\ .\] 
  
Let $\chi$ denote the Euler characteristic for the cohomology groups of $\mathcal{L}$ on the various open subsets, then  $\chi (\mathbb{P}^1)=\chi (X_0)+\chi (Y)-\chi (Y^*)$. 

The restriction of $\mathcal{L}$ to $X_0$ is a locally constant sheaf of $\mathbb{C}$-vector spaces of dimension $m$. This corresponds to a representation  of the fundamental group
$\pi _1(X_0)$ on a vector space $\mathbb{C}^m$. The group $\pi _1(X_0)$ is free on $r-1$
generators $\gamma _1,\dots ,\gamma _{r-1}$. The cohomology groups that we want to calculate
coincide with the group cohomology of the above representation. The latter are the cohomology groups of the complex $0\rightarrow \mathbb{C}^m\rightarrow (\mathbb{C}^m)^{r-1}\rightarrow 0$
where the non trivial map is given by $v\mapsto (\gamma _1v-v,\dots ,\gamma _{r-1}v-v)$.
It follows that $\chi (X_0)=m-(r-1)m=m(2-r)$ and moreover $H^2(X_0,\mathcal{L})=0$.  Using the local calculations we find the formula
 \[\chi (\mathbb{P}^1)= m(2-r)+\sum _{i=1}^rt_0(p_i)\ .  \]
 Suppose that $t_0(p_i)=0$ for some $i$, then $H^i(\mathbb{P}^1,\mathcal{L})=0$
 for $i=0,2$. Thus we proved the following result.
 
 \begin{proposition} Let $L$ be a connection on $\mathbb{P}^1$ of rank $m$ having $r>0$ regular singular points $p_1,\dots ,p_r$ (and no other singularities). Let $t_0(p_i)$
 denote the dimension of the local solution space at the point $p_i$. Suppose that
 at least one $t_0(p_i)$ is zero. Then $lgl(L)=(r-2)m-\sum _{i=1}^rt_0(p_i)$.  
 \end{proposition}

 {\bf References}\\
\noindent [H] G.-J.~van der Heiden -- {\it Weil Pairing and the Drinfeld Modular Curve}
-- Thesis University of Groningen, October 2003, ISBN 90-367-1889-9. \\
\noindent [Jac] N.~Jacobson -- {\it Lie Algebras} -- Interscience Tracts 10, 1962.\\ 
\noindent [Jan] J.C.~Jantzen -- {\it Representations of Algebraic groups} --  second edition -- Mathematical Surveys and Monographs Vol 107, AMS, 2003.\\
\noindent [M] B.~Malgrange -- {\it Sur les points singuliers des \'equations 
diff\'erientielles} -- Ens. Math., 20:149-176, 1974.\\
\noindent [vdP-R]  M. van der Put and M.~Reversat -- {\it Krichever modules for difference and differential equations} --  In: Analyse complexe, syst\`emes  dynamiques, sommabilit\'e des s\'eries divergentes et th\'eories galoisiennes (I). Volume en l'honneur de Jean-Pierre Ramis. Mich\`ele
Loday-Richaud (\'ed.). Ast\'erisque 296, p.207-225, 2004. \\
\noindent [vdP-S] M. van der Put and M.F.~Singer -- {\it Galois theory of linear differential 
equations } -- Grundlehren der mathematische Wissenschaften. Vol 328, Springer Verlag 2003.\\
\noindent [T] P.~Tauvel and Ruppert W.T. Yu - { \it Lie Algebras and Lie Groups} - Springer 
Verlag 2005.\\

 \end{document}